\documentclass[a4paper]{article}
\usepackage{amssymb}
\usepackage{amsmath}

\input{tcilatex}

\begin{document}

\title{On Schubert's Problem of Characteristics}
\author{Haibao Duan and Xuezhi Zhao \thanks{%
The authors are supported by NSFC 11131008, 11431009 and 11661131004.}}
\maketitle

\begin{abstract}
The Schubert varieties on a flag manifold $G/P$ give rise to a cell
decomposition on $G/P$ whose Kronecker duals, known as the Schubert classes
on $G/P$, form an additive base of the integral cohomology $H^{\ast }(G/P)$.
The Schubert's problem of characteristics asks to express a monomial in the
Schubert classes as a linear combination in the Schubert basis.

We present a unified formula expressing the characteristics of a flag
manifold $G/P$ as polynomials in the Cartan numbers of the group $G$. As
application we develop a direct approach to our recent works on the Schubert
presentation of the cohomology rings of flag manifolds $G/P$.

\begin{description}
\item[2010 Mathematical Subject Classification] 14M15; 55T10

\item[Key words and phrases:] Lie group, flag manifold, Schubert variety,
cohomology
\end{description}
\end{abstract}

\section{The problem}

\begin{quote}
\textsl{Schubert considered what he called the problem of characteristics to
be the main theoretical problem of enumerative geometry. }-S. Kleimann \cite[%
1987]{Kl2}

\textsl{The existence of a finite basis for the homologies in every closed
manifold implies furthermore the solvability of Schubert's \textquotedblleft
characteristics problems\textquotedblright\ in general. }-Van der Waerden 
\cite[1930]{Wa}\footnote{%
This English translation of the profound discovery of van der Waerden in
German \cite{Wa} is due to N. Schappacher \cite{NSch}.}
\end{quote}

Schubert calculus is the intersection theory of the 19th century, together
with applications to enumerative geometry. Justifying this calculus was a
major topic of the 20 century algebraic geometry, and was also the content
of Hilbert's 15th problem \textquotedblleft Rigorous foundation of
Schubert's enumerative calculus\textquotedblright\ \cite{Hi,S}. Thanks to
the pioneer works \cite{Wa,E} of Van der Waerden and Ehresmann, the problem
of characteristics \cite{Sc3}, considered by Schubert as the fundamental
problem of enumerative geometry, has had a concise statement by the 1950's.

Let $G$ be compact connected Lie group with a maximal torus $T$. For an one
parameter subgroup $\alpha :\mathbb{R}\rightarrow G$ its centralizer $P$ is
a parabolic subgroup on $G$, while the homogeneous space $G/P$ is a
projective variety, called a \textsl{flag manifold} of $G$. Let $W(P,G)$ be
the set of left cosets of the Weyl group $W$ of $G$ by the Weyl group $W(P)$
of $P$ with associated length function $l:W(P,G)\rightarrow \mathbb{Z}$. The
following result was discovered by Ehresmann \cite[1934]{E} for the
Grassmannians $G_{n,k}$ of $k$-dimensional linear subspaces on $\mathbb{C}%
^{n}$, announced by Chevalley \cite[1958]{Ch} for the complete flag
manifolds $G/T$, and extended to all flag manifolds $G/P$ by
Bernstein-Gel'fand-Gel'fand \cite[1973]{BGG}.

\bigskip

\noindent \textbf{Theorem 1.1. }\textsl{The flag manifold }$G/P$\textsl{\
admits a decomposition into the cells indexed by the elements of }$W(P,G)$%
\textsl{,}

\begin{enumerate}
\item[(1.1)] $G/P=\underset{w\in W(P,G)}{\cup }X_{w}$\textsl{, }$\quad \dim
X_{w}=2l(w)$\textsl{,}
\end{enumerate}

\noindent \textsl{with each cell }$X_{w}$\textsl{\ the closure of an
algebraic affine space, called the Schubert variety on }$G/P$\textsl{\
associated to }$w$\textsl{.}\hfill $\square$

\bigskip

Since only even dimensional cells are involved in the decomposition (1.1),
the set $\{[X_{w}],$ $w\in W(P,G)\}$ of fundamental classes forms an
additive basis of the integral homology $H_{\ast }(G/P)$. The cocycle
classes $s_{w}\in H^{\ast }(G/P)$ Kronecker dual to the basis (i.e. $%
\left\langle s_{w},[X_{u}]\right\rangle =\delta _{w,u}$, $w,u\in W(P,G)$)
gives rise to the \textsl{Schubert class }associated to\textsl{\ }$w\in
W(P,G)$. Theorem 1.1 implies the following result, well known as \textsl{the}
\textsl{basis theorem of Schubert calculus }\cite[\S 8]{Wa}.

\bigskip

\noindent \textbf{Theorem 1.2.} \textsl{The set }$\{s_{w},w\in W(P,G)\}$ 
\textsl{of Schubert classes\ forms an additive basis of the integral
cohomology }$H^{\ast }(G/P)$\textsl{.}\hfill $\square $

\bigskip

An immediate consequence is that any monomial $s_{w_{1}}\cdots s_{w_{k}}$ in
the Schubert classes on $G/P$ can be expressed as a linear combination of
the basis elements

\begin{enumerate}
\item[(1.2)] $s_{w_{1}}\cdots s_{w_{k}}=\sum\limits_{w\in
W(P,G),l(w)=l(w_{1})+\cdots +l(w_{k}),\text{ }}a_{w_{1},\ldots
,w_{k}}^{w}\cdot s_{w}$, $a_{w_{1},\ldots ,w_{k}}^{w}\in \mathbb{Z}$,
\end{enumerate}

\noindent where the coefficients $a_{w_{1},\ldots ,w_{k}}^{w}$ are called 
\textsl{characteristics} by Schubert \cite{Sc3,Wa,E}.

\begin{quote}
\textbf{The problem of characteristics.} Given a monomial $s_{w_{1}}\cdots
s_{w_{k}}$ in the Schubert classes, determine the characteristics $%
a_{w_{1},\ldots ,w_{k}}^{w}$ for all $w\in W(P,G)$ with $l(w)=l(w_{1})+%
\cdots +l(w_{r})$.\hfill $\square $
\end{quote}

The characteristics are of particular importance in geometry, algebra and
topology. They provide solutions to the problems of enumerative geometry 
\cite{Sc1,Sc2,Sc3,Sc4}; were seen by Hilbert as \textquotedblleft the degree
of the final equations and the multiplicity of their
solutions\textquotedblright\ of a system \cite{Hi}; and are requested by
describing the cohomology ring $H^{\ast }(G/P)$ in the Schubert basis \cite[%
p.331]{W}. Notably, the degree of a Schubert variety \cite{Sc4} and the
multiplicative rule of two Schubert classes \cite{Sc3} are two special cases
of the problem which have received considerable attentions in literatures,
see \cite{D0,D} for accounts on the earlier relevant works.

This paper summaries and simplifies our series works \cite{D0,D,DZ1,DZ2,DZ3}
devoted to describe the integral cohomologies of flag manifolds by a minimal
system of generators and relations in the Schubert classes. Precisely, based
on a formula of the characteristics $a_{w_{1},_{\cdots },w_{k}}^{w}$ stated
in Section \S 2 and established in \S 3, we address in Sections \S 4 and \S %
5 a more direct approach to the Schubert presentations \cite{DZ2,DZ3} of the
cohomology rings of flag manifolds $G/P$.

\section{The formula of the characteristics $a_{w_{1},\ldots ,w_{k}}^{w}$}

To investigate the topology of a flag manifold $G/P$ we may assume that the
Lie group $G$ is $1$-connected and simple. Resorting to the geometry of the
Stiefel diagram of the Lie group $G$ we present in Theorem 2.4 a formula
that boils down the characteristics $a_{w_{1},\ldots ,w_{k}}^{w}$ to the
Cartan matrix of the group $G$.

Fix a maximal torus $T$ on $G$ and set $n=\dim T$. Equip the Lie algebra $%
L(G)$ with an inner product $(,)$ so that the adjoint representation acts as
isometries of $L(G)$. The\textsl{\ Cartan subalgebra} of $G$ is the linear
subspace $L(T)\subset L(G)$.

The restriction of the exponential map $\exp :L(G)\rightarrow G$ on $L(T)$
defines a set $S(G)$ of $\frac{1}{2}(\dim G-n)$ hyperplanes on $L(T)$,
namely, the set of\textsl{\ singular hyperplanes }through the origin in $%
L(T) $ \cite[p.226]{BT}. These planes divide $L(T)$ into finitely many
convex regions, each one is called a \textsl{Weyl chamber} of $G$. The
reflections on $L(T)$ in these planes generate the Weyl group $W$ of $G$ 
\cite[p.49]{Hu}.

The map $\exp $ on $L(T)$ carries the normal line $l$ (through the origin)
of a hyperplane $L\in $\ $S(G)$ to a circle subgroup on $T$. Let $\pm \alpha
\in l$ be the non-zero vectors with minimal length so that $\exp (\pm \alpha
)=e$ (the group unit). The set $\Phi (G)$ consisting of all those vectors $%
\pm \alpha $ is called \textsl{the root system of }$G$. Fixing a regular
point $x_{0}\in L(T)-\cup _{{L\in D(G)}}L$ the set of \textsl{simple roots}
relative to $x_{0}$ is

\begin{quote}
$\Delta (x_{0})=\{\beta \in \Phi (G)\mid (\beta ,x_{0})>0\}$ (\cite[p.47]{Hu}%
).
\end{quote}

\noindent In addition, for a simple root $\beta \in \Delta $ \textsl{the} 
\textsl{simple reflection} relative to $\beta $ is the reflection $\sigma _{{%
\beta }}$ in the plane $L_{{\beta }}\in S(G)$ perpendicular to $\beta $. If $%
\beta ,\beta ^{\prime }\in \Delta $ \textsl{the Cartan number }

\begin{quote}
$\beta \circ \beta ^{\prime }:=2(\beta ,\beta ^{\prime })/(\beta ^{\prime
},\beta ^{\prime })$
\end{quote}

\noindent is always an integer, and only $0,\pm 1,\pm 2,\pm 3$ can occur (%
\cite[p.55]{Hu}).

Since the set of simple reflections $\{\sigma _{{\beta }}\mid \beta \in
\Delta \}$ generates $W$ (\cite[p.193]{BT}) every $w\in W$ admits a
factorization of the form

\begin{enumerate}
\item[(2.1)] \noindent $w=\sigma _{{\beta }_{{1}}}\circ \cdots \circ \sigma
_{{\beta }_{{m}}}$, $\beta _{{i}}\in \Delta $.
\end{enumerate}

\noindent \textbf{Definition 2.1}. The \textsl{length} $l(w)$ of an element $%
w\in W$ is the least number of factors in all decompositions of $w$ in the
form (2.1). The decomposition (2.1) is called \textsl{reduced} if $m=l(w)$.

For a reduced decomposition (2.1) of $w$ the $m\times m$ (strictly upper
triangular) matrix $A_{{w}}=(a_{{i,j}})$ with $a_{{i,j}}=0$ if $i\geq j$ and 
$-\beta _{{j}}\circ \beta _{{i}}$ if $i<j$ is called \textsl{the Cartan
matrix of }$w$\textsl{\ }relative to the decomposition\textsl{\ }%
(2.1).\hfill $\square $

\bigskip

\noindent \textbf{Example 2.2.} In \cite{St} Stembridge asked for an
approach to find a reduced decomposition (2.1) for each $w\in W$. Resorting
to the geometry of the Cartan subalgebra $L(T)$ this task can be implemented
by the following method.

Picture $W$ as the $W$-orbit $\{w(x_{0})\in L(T)\mid w\in W\}$ through the
regular point $x_{0}$. For a $w\in W$ let $C_{w}$ be a line segment on $L(T)$
from the Weyl chamber containing $x_{0}$ to $w(x_{0})$, that crosses the
planes in $S(G)$ once at a time. Assume that they are met in the order $%
L_{\alpha _{1}},\ldots ,L_{\alpha _{k}}$, $\alpha _{i}\in \Phi (G)$. Then $%
l(w)=k$ and $w=\sigma _{\alpha _{k}}\circ \cdots \circ \sigma _{\alpha _{1}}$%
. Set

\begin{quote}
$\beta _{1}=\alpha _{1},\quad \beta _{2}=\sigma _{\alpha _{1}}(\alpha
_{2}),\cdots ,\quad \beta _{k}=\sigma _{\alpha _{1}}\circ \cdots \circ
\sigma _{\alpha _{k-1}}(\alpha _{k})$.
\end{quote}

\noindent Then, from $\beta _{i}\in \Delta $ and $\sigma _{\beta
_{i}}=\sigma _{\alpha _{1}}\circ \cdots \circ \sigma _{\alpha _{i-1}}\circ
\sigma _{\alpha _{i}}\circ \sigma _{\alpha _{i-1}}\circ \cdots \circ \sigma
_{\alpha _{1}}$ one sees that $w=\sigma _{\beta _{1}}\circ \cdots \circ
\sigma _{\beta _{k}}$, which is reduced because of $l(w)=k$.\hfill $\square$

\bigskip

Let $\mathbb{Z}[x_{{1}},\ldots ,x_{{m}}]=\oplus _{{n\geq 0}}\mathbb{Z}[x_{{1}%
},\cdots ,x_{{m}}]^{(n)}$ be the ring of integral polynomials in $x_{{1}%
},\ldots ,x_{{m}}$, graded by $\deg x_{{i}}=1$.

\bigskip

\noindent \textbf{Definition 2.3. }For a $m\times m$ strictly upper
triangular integer matrix $A=(a_{{i,j}})$ \textsl{the triangular operator} $%
T_{{A}}$ associated to $A$ is the additive homomorphism $T_{A}:\mathbb{Z}[x_{%
{1}},\ldots ,x_{{m}}]^{(m)}\rightarrow \mathbb{Z}$ defined recursively by
the following elimination rules:

\begin{quote}
i) If $m=1$ (consequently $A=(0)$) then $T_{{A}}(x_{{1}})=1$;

ii) If $h\in \mathbb{Z}[x_{{1}},\ldots ,x_{{m-1}}]^{(m)}$ then $T_{{A}}(h)=0$%
;

iii) For any $h\in \mathbb{Z}[x_{{1}},\ldots ,x_{{m-1}}]^{(m-r)}$ with $%
r\geq 1$,

$\qquad T_{{A}}(h\cdot x_{{m}}^{r})=T_{{A^{\prime }}}(h\cdot (a_{{1,m}}x_{{1}%
}+\cdots +a_{{m-1,m}}x_{{m-1}})^{r-1}),$
\end{quote}

\noindent where $A^{\prime }$ is the $(m-1)\times (m-1)$ strictly upper
triangular matrix obtained from $A$ by deleting both of the $m^{th}$ column
and row.

By additivity, $T_{{A}}$ is defined for every $h\in \mathbb{Z}[x_{{1}%
},\ldots ,x_{{m}}]^{(m)}$ using the unique expansion $h=\sum\limits_{{0\leq
r\leq m}}h_{{r}}\cdot x_{{m}}^{r}$ with $h_{{r}}\in \mathbb{Z}[x_{{1}%
},\ldots ,x_{{m-1}}]^{(m-r)}$.\hfill $\square$

\bigskip

For a parabolic subgroup $P$ of $G$ the set $W(P,G)$ of left cosets of $W$
by $W(P)$ can be identified with the subset of $W$

\begin{quote}
$W(P,G)=\{w\in W\mid l(w)\leq l(ww^{\prime }),w^{\prime }\in W(P)\}$ (by 
\cite[5.1]{BGG}),
\end{quote}

\noindent where $l$ is the length function on $W$. Assume that $w=\sigma _{{%
\beta }_{{1}}}\circ \cdots \circ \sigma _{{\beta }_{{m}}}$, $\beta _{{i}}\in
\Delta $, is a reduced decomposition of an element $w\in W(P,G)$ with
associated Cartan matrix $A_{{w}}=(a_{{i,j}})_{{m\times m}}$. For a
multi-index $I=\{i_{{1}},\ldots ,i_{{k}}\}\subseteq \{1,\ldots ,m\}$ we put $%
\left\vert I\right\vert :=k$ and set

\begin{quote}
$\sigma _{{I}}:=\sigma _{{\beta }_{{i}_{{1}}}}\circ \cdots \circ \sigma _{{%
\beta }_{{i}_{{k}}}}\in W$, $x_{I}:=x_{i_{1}}\cdots x_{i_{k}}\in \mathbb{Z}%
[x_{{1}},\ldots ,x_{{m}}]$.
\end{quote}

\noindent Our promised formula for the characteristics is:

\bigskip

\noindent \textbf{Theorem 2.4.} \textsl{For every monomial }$s_{w_{1}}\cdots
s_{w_{k}}$\textsl{\ in the Schubert classes on }$G/P$ \textsl{with} $%
l(w)=l(w_{1})+\cdots +l(w_{k})$\textsl{\ we have}

\begin{enumerate}
\item[(2.2)] $a_{w_{1},\ldots ,w_{k}}^{w}=T_{A_{{w}}}\left( \underset{%
i=1,\ldots ,k}{\Pi }\left( \underset{\sigma _{{I}}=w_{i},\left\vert
I\right\vert =l(w_{i}),I\subseteq \{1,\ldots ,m\}}{\Sigma }x_{I}\right)
\right) $\textsl{.}
\end{enumerate}

\noindent \textbf{Remarks 2.5.} Formula (2.2) reduces the characteristic $%
a_{w_{1},\ldots ,w_{k}}^{w}$ to a polynomial in the Cartan numbers of the
group $G$, hence applies uniformly to all flag manifolds $G/P$.

If $k=2$ the characteristic $a_{w_{1},w_{2}}^{w}$ is well known as \textsl{a
Littlewood-Richardson coefficient}, and the formula (2.2) has been obtained
by Duan in \cite{D}. In \cite{Wi} M. Willems generalizes the formula of $%
a_{w_{1},w_{2}}^{w}$ to the more general context of flag varieties
associated to Kac-Moody groups, and for the equivariant cohomologies.
Recently, A. Bernstein and E. Richmond \cite{BR} obtained also a formula
expressing $a_{w_{1},w_{2}}^{w}$ in the Cartan numbers of $G$. \hfill $%
\square $

\section{Proof of the characteristics formula (2.2)}

In this paper the homologies and cohomologies are over the ring $\mathbb{Z}$
of integers. If $f:X\rightarrow Y$ is a continuous map between two
topological spaces $X$ and $Y$, $f_{{\ast }}$ (resp. $f^{\ast }$) denote the
homology (resp. cohomology) homomorphism induced by $f$. For an oriented
closed manifold $M$ (resp. a connected projective variety) the notion $%
[M]\in H_{{\dim M}}(M)$ stands for the orientation class. In addition, the
Kronecker pairing between cohomology and homology of a space $X$ is written
as $\left\langle ,\right\rangle :H^{\ast }(X)\times H_{{\ast }%
}(X)\rightarrow \mathbb{Z}$. The proof of Theorem 2.4 makes use of the
celebrated $K$-cycles on $G/T$ constructed by Bott and Samelson in \cite%
{BS,BS1}. We begin by recalling the construction of these cycles, as well as
their basic properties developed in \cite{BS,BS1,D0,D}.

For a simple Lie group $G$ fix a regular point $x_{0}\in L(T)$ and let $%
\Delta $ be the set of simple roots relative to $x_{0}$. For a $\beta \in
\Delta $ let $K_{{\beta }}$ be the centralizer of the subspace $\exp (L_{{%
\beta }})$ on $G$, where $L_{{\beta }}\in S(G)$ is the plane perpendicular
to $\beta $. Then $T\subset K_{{\beta }}$ and the quotient $K_{{\beta }}/T$
is diffeomorphic to the $2$-sphere \cite[p.996]{BS1}.

The $2$-sphere $K_{{\beta }}/T$ carries a natural orientation specified as
follows. The Cartan decomposition of the Lie algebra $L(K_{{\beta }})$
relative to the maximal torus $T\subset K_{{\beta }}$ takes the form $L(K_{{%
\beta }})=L(T)\oplus \vartheta _{{\beta }}$, where $\vartheta _{{\beta }%
}\subset L(G)$ is the \textsl{root space} belonging to the root $\beta $ 
\cite[p.35]{Hu}. Taking a non-zero vector $v\in \vartheta _{{\beta }}$ and
letting $v^{\prime }\in \vartheta _{{\beta }}$ be such that $[v,v^{\prime
}]=\beta $, where $[,]$ is the Lie bracket on $L(G)$, then the ordered base $%
\{v,v\prime \}$ furnishes $\vartheta _{{\beta }}$ with an orientation that
is irrelevant to the choices of $v$. The tangent map of the quotient $\pi _{{%
\beta }}:K_{{\beta }}\rightarrow K_{{\beta }}/T$ at the group unit $e\in K_{{%
\beta }}$ maps the $2$-plane $\vartheta _{{\beta }}$ isomorphically onto the
tangent space to the sphere $K_{{\beta }}/T$ at the point $\pi _{{\beta }%
}(e) $. In this manner the orientation $\{v,v\prime \}$ on $\vartheta _{{%
\beta }}$ furnishes the sphere $K_{{\beta }}/T$ with the orientation $\omega
_{{\beta }}=\{\pi _{{\beta }}(v),\pi _{{\beta }}(v\prime )\}$.

For an ordered sequence $\beta _{{1}},\ldots ,\beta _{{m}}\in \Delta $ of $m$
simple roots (repetitions like $\beta _{{i}}=\beta _{{j}}$ may occur) let $%
K(\beta _{{1}},\ldots ,\beta _{{m}})$ be the product group $K_{{\beta }_{{1}%
}}\times \cdots \times K_{{\beta }_{{m}}}$. With $T\subset K_{{\beta }_{{i}%
}} $ the product $T\times \cdots \times T$ ($m$-copies) acts on $K(\beta _{{1%
}},\ldots ,\beta _{{m}})$ by

\begin{quote}
$(g_{{1}},\ldots ,g_{{m}})(t_{{1}},\ldots ,t_{{m}})=(g_{{1}}t_{{1}},t_{{1}%
}^{-1}g_{{2}}t_{{2}},\ldots ,t_{{m-1}}^{-1}g_{{m}}t_{{m}})$.
\end{quote}

\noindent Let $\Gamma (\beta _{{1}},\ldots ,\beta _{{m}})$ be the base
manifold of this principal action that is oriented by the $\omega _{{\beta }%
_{{i}}}$, $1\leq i\leq m$. The point on $\Gamma (\beta _{{1}},\ldots ,\beta
_{{m}})$ corresponding to the point $(g_{{1}},\ldots ,g_{{m}})\in K(\beta _{{%
1}},\ldots ,\beta _{{m}})$ is called $[g_{{1}},\ldots ,g_{{m}}]$.

The integral cohomology of the oriented manifold $\Gamma (\beta _{{1}%
},\ldots ,\beta _{{m}})$ has been determined by Bott and Samelson in \cite%
{BS}. Let $\varphi _{{i}}:K_{{\beta }_{{i}}}/T\rightarrow \Gamma (\beta _{{1}%
},\ldots ,\beta _{{m}})$ be the embedding induced by the inclusion $K_{{%
\beta }_{{i}}}\rightarrow K(\beta _{{1}},\cdots ,\beta _{{m}})$ onto the $%
i^{th}$ factor group, and put

\begin{quote}
$y_{{i}}=\varphi _{{i}_{\ast }}(\omega _{{\beta }_{{i}}})\in H_{2}(\Gamma
(\beta _{{1}},\ldots ,\beta _{{m}}))$, $1\leq i\leq m$.
\end{quote}

\noindent Form the $m\times m$ strictly upper triangular matrix $A=(a_{{i,j}%
})_{{m\times m}}$ by setting

\begin{quote}
$a_{{i,j}}=0$ if $i\geq j$, but $a_{{i,j}}=-2(\beta _{j},\beta _{i})/(\beta
_{i},\beta _{i}^{\prime })$ if $i<j$.
\end{quote}

\noindent It is easy to see from the construction that the set $\{y_{{1}%
},\ldots ,y_{{m}}\}$\ forms a basis of the second homology group $%
H_{2}(\Gamma (\beta _{{1}},\ldots ,\beta _{{m}}))$.

\bigskip

\noindent \textbf{Lemma 3.1 }(\cite{BS1})\textbf{.} \textsl{Let }$x_{{1}%
},\ldots ,x_{{m}}\in $\textsl{\ }$H^{2}(\Gamma (\beta _{{1}},\ldots ,\beta _{%
{m}}))$ \textsl{be the Kronecker duals of the cycle classes }$y_{{1}},\ldots
,y_{{m}}$\textsl{\ on }$\Gamma (\beta _{{1}},\ldots ,\beta _{{m}})$\textsl{.}
\textsl{Then}

\begin{enumerate}
\item[(3.1)] $H^{\ast }(\Gamma (\beta _{{1}},\ldots ,\beta _{{m}}))=\mathbb{Z%
}[x_{{1}},\ldots ,x_{{m}}]/J$\textsl{,}
\end{enumerate}

\noindent \textsl{where }$J$\textsl{\ is the ideal generated by }$x_{{j}%
}^{2}-\underset{i<j}{\Sigma }a_{{i,j}}x_{{i}}x_{{j}}$\textsl{, }$1\leq j\leq
m$\textsl{.}\hfill $\square $

\bigskip

In view of (3.1) the map $p_{{\Gamma (\beta }_{{1}}{,\ldots ,\beta }_{{m}}{)}%
}$ from the polynomial ring $\mathbb{Z}[x_{{1}},\ldots ,x_{{m}}]$ onto its
quotient $H^{\ast }(\Gamma (\beta _{{1}},\ldots ,\beta _{{m}}))$ gives rise
to the additive map

\begin{quote}
$\int_{{\Gamma (\beta }_{{1}}{,\ldots ,\beta }_{{m}}{)}}:\mathbb{Z}[x_{{1}%
},\ldots ,x_{{m}}]^{(m)}\rightarrow \mathbb{Z}$
\end{quote}

\noindent evaluated by $\int_{{\Gamma (\beta }_{{1}}{,\cdots ,\beta }_{{m}}{)%
}}h=\left\langle p_{{\Gamma (\beta }_{{1}}{,\ldots ,\beta }_{{m}}{)}%
}(h),[\Gamma (\beta _{{1}},\ldots ,\beta _{{m}})]\right\rangle $. The
geometric implication of the triangular operator $T_{{A}}$ in Definition 2.3
is shown by the following result.

\bigskip

\noindent \textbf{Lemma 3.2 }(\cite[Proposition 2]{D0}). \textsl{We have}

\begin{quote}
$\int_{{\Gamma (\beta }_{{1}}{,\ldots ,\beta }_{{m}}{)}}=T_{{A}}:\mathbb{Z}%
[x_{{1}},\ldots ,x_{{m}}]^{(m)}\rightarrow \mathbb{Z}$\textsl{.\ }
\end{quote}

\noindent \textsl{In particular, }$\int_{{\Gamma (\beta }_{{1}}{,\ldots
,\beta }_{{m}}{)}}x_{{1}}\cdots x_{{m}}=1$\textsl{.}\hfill $\square$

\bigskip

For a parabolic $P$ on $G$ we can assume, without loss of the generalities,
that $T\subseteq P\subset G$. For a sequence $\beta _{{1}},\ldots ,\beta _{{m%
}}$ of simple roots \textsl{the associated Bott-Samelson's} \textsl{K-cycle}
on $G/P$ is the map

\begin{quote}
$\varphi _{{\beta }_{{1}}{,\ldots ,\beta }_{{m}}{;P}}:\Gamma (\beta _{{1}%
},\ldots ,\beta _{{m}})\rightarrow G/P$
\end{quote}

\noindent defined by $\varphi _{{\beta }_{{1}}{,\ldots ,\beta }_{{m}}{;P}%
}([g_{{1}},\ldots ,g_{{m}}])=g_{{1}}\cdots g_{{m}}P$. If $P=T$ Hansen \cite%
{H} has shown that certain $K$-cycles are desingularizations of the Schubert
varieties on $G/T$. The following more general result allows one to
translate the calculation with Schubert classes on $G/P$ to computing with
monomials in the much simpler ring $H^{\ast }(\Gamma (\beta _{{1}},\ldots
,\beta _{{m}}))$.

\bigskip

\noindent \textbf{Lemma 3.3.} \textsl{With respect the Schubert basis on} $%
H^{\ast }(G/P)$ \textsl{the induced map of }$\varphi _{{\beta }_{{1}}{%
,\ldots ,\beta }_{{m}}{;P}}$\textsl{\ on the cohomologies is given by}

\begin{enumerate}
\item[(3.2)] $\varphi _{\beta _{1},\ldots ,\beta _{m};P}^{\ast
}(s_{w})=(-1)^{l(w)}\underset{_{\sigma _{I}=w,\left\vert I\right\vert
=l(w),I\subseteq \lbrack 1,\ldots ,m]}}{\Sigma }x_{I}$\textsl{,} $w\in
W(P;G) $\textsl{.}
\end{enumerate}

\noindent \textbf{Proof. }With $T\subseteq P\subset G$ the map $\varphi
_{\beta _{1},\ldots ,\beta _{m};P}$ factors through $\varphi _{\beta
_{1},\ldots ,\beta _{m};T}$ in the fashion

\begin{quote}
$%
\begin{array}{ccc}
\Gamma (\beta _{1},\ldots ,\beta _{m}) & \overset{\varphi _{\beta
_{1},\ldots ,\beta _{m};T}}{\rightarrow } & G/T \\ 
& \searrow & \downarrow \pi \\ 
\varphi _{\beta _{1},\ldots ,\beta _{m};P} &  & G/P%
\end{array}%
$,
\end{quote}

\noindent where the map $\pi $ is the fibration with fiber $P/T$. By \cite[%
Lemma 5.1]{D}\textsl{\ }formula (3.2) holds for the case $P=T$. According to 
\cite[\S 5]{BGG}\textsl{\ }the induced map $\pi ^{\ast }:H^{\ast
}(G/P)\rightarrow H^{\ast }(G/T)$\ is given by $\pi ^{\ast }(s_{w})=s_{w}$,%
\textsl{\ }$w\in W(P;G)$,\textsl{\ }showing formula (3.2) for the general
case $T\subset P$.\hfill $\square $

\bigskip

\noindent \textbf{Proof of Theorem 2.4. }For a monomial $s_{w_{1}}\cdots
s_{w_{k}}$ in the Schubert classes of $G/P$ assume as in (1.2) that

\begin{enumerate}
\item[(3.3)] $s_{w_{1}}\cdots s_{w_{k}}=\sum\limits_{w\in
W(P;G),l(w)=m}a_{w_{1},\ldots ,w_{k}}^{w}\cdot s_{w}$, $a_{w_{1},\ldots
,w_{k}}^{w}\in \mathbb{Z}$,
\end{enumerate}

\noindent where $m=l(w_{1})+\cdots +l(w_{k})$. For an element $w_{0}\in
W(P;G)$ with a reduced decomposition $w_{0}=\sigma _{{\beta }_{{1}}}\circ
\cdots \circ \sigma _{{\beta }_{{m}}}$, $\beta _{{i}}\in \Delta $, let $A_{{w%
}_{0}}=(a_{{i,j}})_{{m\times m}}$ be the relative Cartan matrix. Applying
the ring map $\varphi _{\beta _{1},\ldots ,\beta _{m};P}^{\ast }$ to the
equation (3.3) on $H^{\ast }(G/P)$ we obtain by (3.2) the equality on the
group $H^{2m}(\Gamma (\beta _{1},\ldots ,\beta _{m}))$

\begin{center}
$(-1)^{l(w_{1})+\cdots +l(w_{k})}\underset{1\leq i\leq k}{\Pi }(\underset{%
_{\sigma _{I}=w_{i},\left\vert I\right\vert =l(w_{i}),I\subseteq \left\{
1,\ldots ,m\right\} }}{\Sigma }x_{I})=(-1)^{m}a_{w_{1},\ldots
,w_{k}}^{w_{0}}\cdot x_{{1}}\cdots x_{{m}}$.
\end{center}

\noindent Applying $\int_{{\Gamma (\beta }_{{1}}{,\ldots ,\beta }_{{m}}{)}}$
to both sides we get by Lemma 3.2\ that

\begin{center}
$(-1)^{l(w_{1})+\cdots +l(w_{k})}\cdot T_{A_{w}}(\underset{1\leq i\leq k}{%
\Pi }(\underset{_{\sigma _{I}=w_{i},\left\vert I\right\vert
=l(w_{i}),I\subseteq \left\{ 1,\ldots ,m\right\} }}{\Sigma }%
x_{I}))=(-1)^{m}\cdot a_{w_{1},\ldots ,w_{k}}^{w_{0}}$.
\end{center}

\noindent This is identical to (2.2) because of $m=l(w_{1})+\cdots +l(w_{k})$%
.\hfill $\square$

\section{The cohomology of flag manifolds $G/P$}

\begin{quote}
\textsl{The classical Schubert calculus amounts to the determination of the
intersection rings on Grassmann varieties and on the so called "flag
manifolds" of projective geometry. }-A. Weil \cite[p.331]{W}
\end{quote}

A classical problem of topology is to express the integral cohomology ring $%
H^{\ast }(G/H)$ of a homogeneous space $G/H$ by a minimal system of explicit
generators and relations. The traditional approach due to H. Cartan, A.
Borel, P. Baum, H. Toda utilize various spectral sequence techniques \cite%
{B,B1,HMS,T,Wo}, and the calculation encounters the same difficulties when
applied to a Lie group $G$ with torsion elements in its integral cohomology,
in particular, when $G$ is one of the five exceptional Lie groups \cite%
{IT,T,W2}.

However, if $P\subset G$ is parabolic, Schubert calculus makes the structure
of the ring $H^{\ast }(G/P)$ appearing in a new light. Given a set $\left\{
y_{{1}},\ldots ,y_{{k}}\right\} $ of $k$ elements let $\mathbb{Z}[y_{{1}%
},\ldots ,y_{{k}}]$ be the ring of polynomials in $y_{{1}},\ldots ,y_{{k}}$
with integer coefficients. For a subset $\left\{ r_{1},\ldots ,r_{m}\right\} 
$ $\subset \mathbb{Z}[y_{{1}},\ldots ,y_{{k}}]$ of homogenous polynomials
denote by $\langle r_{1},\ldots ,r_{m}\rangle $ the ideal generated by $%
r_{1},\ldots ,r_{m}$.

\bigskip

\noindent \textbf{Theorem 4.1.} \textsl{For each flag manifold }$G/P$\textsl{%
\ there exist a set }$\left\{ y_{{1}},\ldots ,y_{{k}}\right\} $ \textsl{of
Schubert classes on} $G/P$\textsl{, and a set }$\left\{ r_{1},\ldots
,r_{m}\right\} $\textsl{\ }$\subset \mathbb{Z}[y_{{1}},\ldots ,y_{{k}}]$%
\textsl{\ of polynomials,} \textsl{so that the inclusion }$\left\{ y_{{1}%
},\ldots ,y_{{k}}\right\} \subset H^{\ast }(G/P)$\textsl{\ induces a ring
isomorphism}

\begin{enumerate}
\item[(4.1)] $H^{\ast }(G/P)=\mathbb{Z}[y_{{1}},\ldots ,y_{{k}}]/\langle r_{{%
1}},\ldots ,r_{{m}}\rangle $\textsl{,}
\end{enumerate}

\noindent \textsl{where both the numbers }$k$\textsl{\ and }$m$\textsl{\ are
minimal subject to this presentation.}\hfill $\square$

\bigskip

\noindent \textbf{Proof.} Let $D(H^{\ast }(G/P))\subset H^{\ast }(G/P)$ be
the ideal of the decomposable elements. Since the ring $H^{\ast }(G/P)$ is
torsion free and has a basis consisting of Schubert classes, there is a set $%
\left\{ y_{{1}},\ldots ,y_{{k}}\right\} $ of Schubert classes on $G/P$ that
corresponds to a basis of the quotient group $H^{\ast }(G/P)/D(H^{\ast
}(G/P))$. In particular, the inclusion $\left\{ y_{{1}},\ldots ,y_{{k}%
}\right\} \subset H^{\ast }(G/P)$ induces a surjective ring map

\begin{quote}
$f:\mathbb{Z}[y_{{1}},\ldots ,y_{{k}}]\rightarrow H^{\ast }(G/P)$.
\end{quote}

\noindent Since $\ker f$ is an ideal the Hilbert basis theorem implies that
there exists a finite subset $\{r_{{1}},\ldots ,r_{{m}}\}\subset \mathbb{Z}%
[y_{{1}},\ldots ,y_{{n}}]$ so that $\ker f=$\textsl{\ }$\left\langle r_{{1}%
},\ldots ,r_{{m}}\right\rangle $. We can of course assume that the number $m$
is minimal subject to this constraint.

As the cardinality of a basis of the quotient group $H^{\ast
}(G/P)/D(H^{\ast }(G/P))$ the number $k$ is an invariant of $G/P$. In
addition, if one changes the generators $y_{{1}},\ldots ,y_{{k}}$ to $y_{{1}%
}^{\prime },\ldots ,y_{{k}}^{\prime }$, then each old generator $y_{{i}}$
can be expressed as a polynomial $g_{i}$ in the new ones $y_{{1}}^{\prime
},\ldots ,y_{{k}}^{\prime }$, and the invariance of the number $m$ is shown
by the presentation

\begin{quote}
$H^{\ast }(G/P)=\mathbb{Z}[y_{{1}}^{\prime },\ldots ,y_{{k}}^{\prime
}]/\langle r_{1}^{\prime },\ldots ,r_{m}^{\prime }\rangle $,
\end{quote}

\noindent where $r_{j}^{\prime }$ is obtained from $r_{j}$ by substituting $%
g_{i}$ for $y_{{i}}$, $1\leq j\leq m$.\hfill $\square$

\bigskip

A presentation of the ring $H^{\ast }(G/P)$ in the form of (4.1) will be
called a \textsl{Schubert presentation} \textsl{of the cohomology of} $G/P$,
while the set $\left\{ y_{{1}},\ldots ,y_{{k}}\right\} $ of generators will
be called \textsl{a set of special Schubert classes} on $G/P$. Based on the
characteristic formula (2.2) we develop in this section algebraic and
computational machineries implementing Schubert presentation of the ring $%
H^{\ast }(G/P)$. To be precise the following conventions will be adopted
throughout the remaining part of this section.

\begin{quote}
i) $G$ is a $1$-connected simple Lie group with Weyl group $W$, and a fixed
maximal torus $T$;

ii) A set $\Delta =\{\beta _{1},\ldots ,\beta _{n}\}$ of simple roots of $G$
is given and ordered as the vertex of the Dykin diagram of $G$ pictured on 
\cite[p.58]{H};

iii) For each simple root $\beta _{i}\in \Delta $ write $\sigma _{i}$
instead of $\sigma _{\beta _{i}}\in W$; $\omega _{i}$ in place of the
Schubert class $s_{\sigma _{\beta _{i}}}\in H^{2}(G/T)$.
\end{quote}

\noindent Note that Theorem 1.2 implies that the set $\{\omega _{1},\ldots
,\omega _{n}\}$ is the Schubert basis of the second cohomology $H^{2}(G/T)$,
whose elements is identical to \textsl{the fundamental} \textsl{dominant
weights} of $G$ in the context of Borel and Hirzebruch \cite{BH,Du}.

\subsection{Decomposition}

By convention iii) each $w\in W$ admits a factorization of the form

\begin{quote}
$w=\sigma _{{i}_{{1}}}\circ \cdots \circ \sigma _{{i}_{{k}}}$, $1\leq i_{{1}%
},\ldots ,i_{{k}}\leq n$, $l(w)=k$.
\end{quote}

\noindent hence can be written as $w=\sigma \lbrack I]$ with $I=(i_{{1}%
},\ldots ,i_{{k}})$. Such expressions of $w$ may not be unique, but the
ambiguity can be dispelled by employing the following notion. Furnish the
set of all reduced decompositions of $w$

\begin{quote}
$D(w):=\{w=\sigma \lbrack I]\mid I=(i_{{1}},\ldots ,i_{{k}}),l(w)=k\}$.
\end{quote}

\noindent with the order $\leq $ given by the lexicographical order on the
multi-indexes $I$. Call a decomposition $w=\sigma \lbrack I]$ \textsl{%
minimized }if $I$ is the minimal one with respect to the order. As result
every $w\in W$\ possesses a unique minimized decomposition.

For a subset $K\subset \{1,\ldots ,n\}$ let $P_{{K}}\subset G$ be the
centralizer of the $1$-parameter subgroup $\{\exp (tb)\in G\mid t\in \mathbb{%
R}\}$ on $G$, where $b\in L(T)$ is a vector that satisfies

\begin{quote}
$\left( \beta _{i},b\right) >0$ if $i\in K$; $\left( \beta _{i},b\right) =0$
if $i\notin K$.
\end{quote}

\noindent Then every parabolic subgroup $P$ is conjugate in $G$ to some $P_{{%
K}}$ with $K\subset \{1,\ldots ,n\}$, while the Weyl group $W(P)\subset W$
is generated by the simple reflections $\sigma _{{j}}$ with $j\notin K$.
Resorting to the length function $l$ on $W$ we embed the set $W(P;G)$ as the
subset of $W$ (as in Section 2)

\begin{quote}
$W(P;G)=\{w\in W\mid l(w_{{1}})\geq l(w)$, $w_{{1}}\in wW(P)\}$,
\end{quote}

\noindent and put $W^{r}(P;G):=\{w\in W(P;G)\mid l(w)=r\}$. Since every $%
w\in W^{r}(P;G)$ has a unique minimized decomposition as $w=\sigma \lbrack
I] $, the set $W^{r}(P;G)$ is also ordered by the lexicographical order on
the multi-index $I$'s, hence can be expressed as

\begin{enumerate}
\item[(4.2)] $W^{r}(P;G)=\{w_{{r,i}}\mid 1\leq i\leq \beta (r)\}$, $\beta
(r):=\left\vert W^{r}(P;G)\right\vert $,
\end{enumerate}

\noindent where $w_{{r,i}}$ is the $i^{th}$ element in the ordered set $%
W^{r}(P;G)$. In \cite{DZ1} a program entitled \textquotedblleft \textsl{%
Decomposition}\textquotedblright\ is composed, whose function is summarized
below.

\begin{quote}
\noindent \textbf{Algorithm 4.2.} \textsl{Decomposition.}

\textbf{Input:} \textsl{The Cartan matrix }$A=(a_{{ij}})_{{n\times n}}$ 
\textsl{of }$G$\textsl{, and a subset }$K\subset \{1,\ldots ,n\}$\textsl{.}

\textbf{Output: }\textsl{The set }$W(P_{K};G)$\textsl{\ being presented by
the minimized decompositions of its elements, together with the index system
(4.2) imposed by the order }$\leq $ \textsl{.}
\end{quote}

\noindent For examples of the results coming from \textsl{Decomposition} we
refer to \cite[1.1--7.1]{DZ4}.

\subsection{Factorization of the ring $H^{\ast }(G/T)$ using fibration}

The cardinality of the Schubert basis of $G/T$ agrees with the order of the
Weyl group $W$, which in general is very large. To reduce the computation
costs we may take a proper subset $K\subset \{1,\ldots ,n\}$ and let $%
P:=P_{K}$ be the corresponding parabolic subgroup. The inclusion $T\subset P$
$\subset G$ then induces the fibration

\begin{enumerate}
\item[(4.3)] $P/T\overset{i}{\hookrightarrow }G/T\overset{\pi }{\rightarrow }%
G/P$,
\end{enumerate}

\noindent where the induced maps $\pi ^{\ast }$ and $i^{\ast }$ behave well
with respect to the Schubert bases of the three flag manifolds $P/T$, $G/P$
and $G/T$ in the following sense:

\begin{quote}
i) With respect to the inclusion $W(P)\subset W$\ the map $i^{\ast }$\
carries the subset $\{s_{w}\}_{w\in W(P)\subset W}$\ of the Schubert basis
of $H^{\ast }(G/T)$\ onto the Schubert basis of $H^{\ast }(P/T)$.

ii) With respect to the inclusion $W(P;G)\subset W$\ the map $\pi ^{\ast }$
identifies the Schubert basis $\{s_{w}\}_{w\in W(P;G)}$\ of $H^{\ast }(G/P)$%
\ with a subset\ of the Schubert basis of $H^{\ast }(G/T)$.
\end{quote}

\noindent For these reasons we can make no difference in notation between an
element in $H^{\ast }(G/P)$\ and its $\pi ^{\ast }$\ image in $H^{\ast
}(G/T) $, and between a Schubert class on $P/T$\ and its $i^{\ast }$\
pre-image on $G/T$.

Assume now that $\{y_{1},\ldots ,y_{n_{1}}\}$\ and $\{x_{1},\ldots
,x_{n_{2}}\}$ are respectively special Schubert classes on $P/T$ and $G/P$,
and with respect to them one has the Schubert presentations

\begin{enumerate}
\item[(4.4)] $H^{\ast }(P/T)=\frac{\mathbb{Z}[y_{i}]_{1\leq i\leq n_{1}}}{%
\left\langle h_{s}\right\rangle _{1\leq s\leq m_{1}}}$; $H^{\ast }(G/P)=%
\frac{\mathbb{Z}[x_{j}]_{1\leq j\leq n_{2}}}{\left\langle r_{t}\right\rangle
_{1\leq t\leq m_{2}}}$,
\end{enumerate}

\noindent where $h_{s}\in \mathbb{Z}[y_{i}]_{1\leq i\leq n_{1}}$, $r_{t}\in 
\mathbb{Z}[x_{j}]_{1\leq j\leq n_{2}}$. The following result allows one to
formulate the ring $H^{\ast }(G/T)$ by the simpler ones $H^{\ast }(P/T)$ and 
$H^{\ast }(G/P)$.

\bigskip

\noindent \textbf{Theorem 4.3.} \textsl{The inclusions }$y_{i},x_{j}\in
H^{\ast }(G/T)$ \textsl{induces a surjective ring map}

\begin{quote}
\textsl{\ }$\varphi :\mathbb{Z}[y_{i},x_{j}]_{1\leq i\leq n_{1},1\leq j\leq
n_{2}}\rightarrow H^{\ast }(G/T)$\textsl{. }
\end{quote}

\noindent \textsl{Furthermore, if }$\{\rho _{s}\}_{1\leq s\leq m_{1}}\subset 
\mathbb{Z}[y_{i},x_{j}]$\textsl{\ is a system satisfying}

\begin{enumerate}
\item[(4.5)] \textsl{\ }$\rho _{s}\in \ker \varphi $\textsl{\ and }$\rho
_{s}\mid _{x_{j}=0}=h_{s}$\textsl{,}
\end{enumerate}

\noindent \textsl{then }$\varphi $\textsl{\ induces a ring isomorphism}

\begin{enumerate}
\item[(4.6)] $H^{\ast }(G/T)=\mathbb{Z}[y_{i},x_{i}]_{1\leq i\leq
n_{1},1\leq j\leq n_{2}}/\left\langle \rho _{s},r_{t}\right\rangle _{1\leq
s\leq m_{1},1\leq t\leq m_{2}}$\textsl{.}
\end{enumerate}

\noindent \textbf{Proof.} By the property i) above the bundle (4.3) has the
Leray-Hirsch property. That is, the cohomology $H^{\ast }(G/T)$ is a free
module over the ring $H^{\ast }(G/P)$ with the basis $\{1,s_{w}\}_{w\in
W(P)} $:

\begin{enumerate}
\item[(4.7)] $H^{\ast }(G/T)=H^{\ast }(G/P)\{1,s_{w}\}_{w\in W(P)}$ (\cite[%
p.231]{Hus}),
\end{enumerate}

\noindent implying that $\varphi $ surjects. It remains to show that for any 
$g\in \ker \varphi $ one has

\begin{quote}
$g\in \left\langle \rho _{s},r_{t}\right\rangle _{1\leq s\leq m_{1},1\leq
t\leq m_{2}}$.
\end{quote}

\noindent To this end we notice by (4.5) and (4.7) that

\begin{quote}
$g\equiv \sum\nolimits_{w\in W(P_{K})}g_{w}\cdot s_{w}\func{mod}
\left\langle \rho _{s}\right\rangle _{1\leq s\leq m_{1}}$ with $g_{w}\in 
\mathbb{Z}[x_{j}]_{1\leq j\leq n_{2}}$.
\end{quote}

\noindent Thus $\varphi (g)=0$ implies $\varphi (g_{w})=0$, showing $%
g_{w}\in \left\langle r_{t}\right\rangle _{1\leq t\leq m_{2}}$ by
(4.4).\hfill $\square $

\subsection{The generalized Grassmannians}

For a topological space $X$ we set

\begin{quote}
$H^{\text{even}}(X):=\oplus _{{r\geq 0}}H^{2r}(X)$,\quad\ $H^{\text{odd}%
}(X):=\oplus _{{r\geq 0}}H^{2r+1}(X)$.
\end{quote}

\noindent Then $H^{\text{even}}(X)$ is a subring of $H^{\ast }(X)$, while $%
H^{\text{odd}}(X)$ is a module over the ring $H^{\text{even}}(X)$.

If $P$ is a parabolic subgroup that corresponds to a singleton $K=\{i\}$,
the flag manifold $G/P$ is called \textsl{generalized Grassmannians of }$G$%
\textsl{\ corresponding to the weight }$\omega _{i}$ \cite{DZ2}. With $%
W^{1}(P;G)=\{\sigma _{i}\}$ consisting of a single element the basis theorem
implies that $H^{2}(G/P)=\mathbb{Z}$\textsl{\ }is generated by\textsl{\ }$%
\omega _{i}$. Furthermore, letting $P^{s}$ be the semi-simple part of $P$,
then the projection $p:G/P^{s}\rightarrow G/P$ is an oriented circle bundle
on $G/P$ with Euler class $\omega _{i}$. With $H^{\text{odd}}(G/P)=0$ by the
basis theorem the Gysin sequence \cite[p.143]{MS}

\begin{quote}
$\cdots \rightarrow H^{r}(G/P)\overset{{\small p}^{\ast }}{{\rightarrow }}%
H^{r}(G/P^{s})\overset{\beta }{{\rightarrow }}H^{r-1}(G/P)\overset{\omega {%
\cup }}{{\rightarrow }}H^{r+1}(G/P)\overset{{\small p}^{\ast }}{{\rightarrow 
}}\cdots $.
\end{quote}

\noindent of $p$ breaks into the short exact sequences

\begin{enumerate}
\item[(4.8)] $0\rightarrow \omega _{i}\cup H^{2r-2}(G/P)\rightarrow
H^{2r}(G/P)\overset{p^{\ast }}{\rightarrow }H^{2r}(G/P^{s})\rightarrow 0$
\end{enumerate}

\noindent as well as the isomorphisms

\begin{enumerate}
\item[(4.9)] $\beta :H^{2r-1}(G/P^{s})\cong \ker \{H^{2r-2}(G/P)\overset{%
\omega _{i}\cup }{\rightarrow }H^{2r}(G/P)\}$,
\end{enumerate}

\noindent where $\omega _{i}\cup $ means taking cup product with $\omega
_{i} $. In particular, formula (4.8) implies that

\textbf{\ }

\noindent \textbf{Lemma 4.4.} \textsl{If }$S=\{y_{{1}},\ldots ,y_{{m}%
}\}\subset H^{\ast }(G/P)$\textsl{\ is a subset so that }$p^{\ast
}S=\{p^{\ast }(y_{{1}}),\ldots ,p^{\ast }(y_{{m}})\}$\textsl{\ is a minimal
set of generators of the ring }$H^{\text{even}}(G/P^{s})$\textsl{, then }$%
S^{\prime }=\{\omega _{i},y_{{1}},\ldots ,y_{{m}}\}$\textsl{\ is a minimal
set of generators of }$H^{\ast }(G/P)$.\hfill $\square$

\bigskip

By Lemma 4.4 the inclusions $\{\omega _{i}\}\cup S$ $\subset H^{\ast }(G/P)$%
, $p^{\ast }S\subset H^{\ast }(G/P^{s})$ extend to the surjective maps $\pi $
and $\overline{\pi }$ that fit in the commutative diagram

\begin{enumerate}
\item[(4.10)] 
\begin{tabular}{llllll}
&  & $\mathbb{Z}[\omega _{i},y_{{1}},\ldots ,y_{{m}}]^{(2r)}$ & $\overset{%
\varphi }{\rightarrow }$ & $\mathbb{Z}[y_{{1}},\ldots ,y_{{m}}]^{(2r)}$ & 
\\ 
&  & $\ \ \ \pi \downarrow $ &  & $\ \ \ \overline{\pi }\downarrow $ &  \\ 
$H^{2r-2}(G/P)$ & $\overset{\omega _{i}\cup }{\rightarrow }$ & $\
H^{2r}(G/P) $ & $\overset{p^{\ast }}{\rightarrow }$ & $\ H^{2r}(G/P^{s})$ & $%
\rightarrow 0$%
\end{tabular}
\end{enumerate}

\noindent where $\mathbb{Z}[\omega _{i},y_{{1}},\ldots ,y_{{m}}]$ is graded
by $\deg \omega _{i}=2,\deg y_{{i}}$, and where

\begin{quote}
$\varphi (\omega _{i})=0$, $\varphi (y_{{i}})=y_{{i}}$; $\overline{\pi }(y_{{%
i}})=p^{\ast }(y_{{i}})$.
\end{quote}

\noindent The next result showing in \cite[Lemma 8]{DZ2} enables us to
formulate a presentation of the ring $H^{\ast }(G/P)$ in term of $H^{\ast
}(G/P^{s})$.

\bigskip

\noindent \textbf{Theorem 4.5. }\textsl{Assume that }$\{h_{{1}},\ldots ,h_{{d%
}}\}\subset \mathbb{Z}[y_{{1}},\ldots ,y_{{m}}]$\textsl{\ is a subset so that%
}

\begin{enumerate}
\item[(4.11)] $H^{\text{even}}(G/P^{s})=\mathbb{Z}[p^{\ast }(y_{{1}}),\ldots
,p^{\ast }(y_{{m}})]/\left\langle p^{\ast }(h_{{1}}),\ldots ,p^{\ast }(h_{{d}%
})\right\rangle $\textsl{,}
\end{enumerate}

\noindent \textsl{and} \textsl{that }$\{d_{{1}},\ldots ,d_{{t}}\}$ \textsl{%
is a basis of the module }$H^{\text{odd}}(G/P^{s})$ \textsl{over }$H^{\text{%
even}}(G/P^{s})$\textsl{. Then}

\begin{enumerate}
\item[(4.12)] $H^{\ast }(G/P)=\mathbb{Z}[\omega _{i},y_{{1}},\ldots ,y_{{m}%
}]/\left\langle r_{{1}},\ldots ,r_{{d}};\omega _{i}g_{{1}},\ldots ,\omega
_{i}g_{{t}}\right\rangle $,\textsl{\ }
\end{enumerate}

\noindent \textsl{where }$\{r_{{1}},\ldots ,r_{{d}}\}$\textsl{, }$\{g_{{1}%
},\ldots ,g_{{t}}\}\subset \mathbb{Z}[\omega _{i},y_{{1}},\ldots ,y_{{m}}]$%
\textsl{\ are two sets of polynomials that satisfy respectively the
following \textquotedblleft \textsl{initial constraints}\textquotedblright }

\begin{quote}
\textsl{i)\ }$r_{{k}}\in \ker \pi $ \textsl{with} $r_{{k}}\mid _{{\omega }%
_{i}{=0}}$ $=h,$\textsl{\quad }$1\leq k\leq d$\textsl{;}

\textsl{ii)\ }$\pi (g_{{j}})=\beta (d_{{j}})$\textsl{, }$1\leq j\leq t$%
\textsl{.}\hfill $\square$
\end{quote}

\subsection{The Characteristics}

To make Theorems 4.3 and 4.5 applicable in practical computation we develop
in this section a series of three algorithms, entitled \textsl{%
Characteristics, Null-space, Giambelli polynomials}, all of them are based
on the characteristic formula (2.2).

\bigskip

\textbf{4.4.1. The Characteristics.}\textsl{\ }For a $w\in W(P;G)$ with the
minimized decomposition $w=\sigma _{{i}_{{1}}}\circ \cdots \circ \sigma _{{i}%
_{{m}}}$, $1\leq i_{{1}},\ldots ,i_{{m}}\leq n$, $l(w)=m$, we observe in
formula (2.2) that

\begin{quote}
i) The Cartan matrix $A_{{w}}$ of $w$ can be read directly from the Cartan
matrix \cite[p.59]{Hu} of the Lie group $G$;

ii) For a $u\in W(P;G)$ with $l(u)=r<m$ the solutions in the multi-indices $%
I=\{j_{{1}},\ldots ,j_{{r}}\}\subseteq \{i_{{1}},\ldots ,i_{{m}}\}$ to the
equation $\sigma _{{I}}=u$ in $W(P;G)$ agree with the solutions to the
linear system $\sigma _{I}(x_{0})=u(x_{0})$ on the vector space $L(T)$,
where $x_{0}$ is the fixed regular point;

iii) The evaluation the operator $T_{{A}_{{w}}}$ on a polynomial have been
programmed using different methods in \cite{DZ1,ZZDL}.
\end{quote}

\noindent Summarizing, granted with \textsl{Decomposition}, formula (2.2)
indicates an effective algorithm to implement a parallel program whose
function is briefed below.

\begin{quote}
\textbf{Algorithm 4.6: }\textsl{Characteristics.}

\textbf{Input:} \textsl{The Cartan matrix }$A=(a_{{ij}})_{{n\times n}}$ 
\textsl{of }$G$\textsl{, and} \textsl{a monomial }$s_{w_{1}}\cdots s_{w_{k}}$%
\textsl{\ in Schubert classes on }$G/P$\textsl{. }

\textbf{Output: }\textsl{The expansion (1.2) of }$s_{w_{1}}\cdots s_{w_{k}}$%
\textsl{\ in the Schubert basis}$.\square $
\end{quote}

\bigskip

\textbf{4.4.2. The Null-space. }Let $\mathbb{Z}[y_{{1}},\ldots ,y_{{k}}]$ be
the ring of polynomials in $y_{{1}},\ldots ,y_{{k}}$ graded by $\deg y_{{i}%
}>0$, and let $\mathbb{Z}[y_{{1}},\ldots ,y_{{k}}]^{(m)}$ be the $\mathbb{Z}$%
-module consisting of all the homogeneous polynomials with degree $m$.
Denote by $\mathbb{N}^{k}$ the set of all $k$-tuples $\alpha =(b_{{1}%
},\ldots ,b_{{k}})$ of non-negative integers. Then \textsl{the set of
monomials basis} of $\mathbb{Z}[y_{{1}},\ldots ,y_{{k}}]^{(m)}$ is

\begin{enumerate}
\item[(4.13)] $B(m)=\{y^{\alpha }=y_{{1}}^{b_{1}}\cdots y_{{k}}^{b_{k}}\mid
\alpha =(b_{{1}},\ldots ,b_{{k}})\in \mathbb{N}^{k},$ $\deg y^{\alpha }=m\}$.
\end{enumerate}

\noindent It will be considered as an ordered set with respect to the
lexicographical order on $\mathbb{N}^{k}$, whose cardinality is called $b(m)$%
.

Let $S=\{y_{{1}},\ldots ,y_{{k}}\}$ be a set of Schubert classes on $G/P$
that generates the ring $H^{\ast }(G/P)$ multiplicatively. Then the
inclusion $S\subset H^{\ast }(G/P)$ induces a surjective ring map $f:\mathbb{%
Z}[y_{{1}},\ldots ,y_{{k}}]\rightarrow H^{\ast }(G/P)$ whose restriction to
degree $2m$ is

\begin{quote}
$f_{{m}}:\mathbb{Z}[y_{{1}},\ldots ,y_{{k}}]^{(2m)}\rightarrow H^{2m}(G/P)$.
\end{quote}

\noindent Combining the \textsl{Characteristics }with the function
\textquotedblleft \textsl{Null-space}\textquotedblright\ in \textsl{%
Mathematica}, a basis of $\ker f_{{m}}$ can be explicitly exhibited.

Let $s_{{m,i}}$ be the Schubert class corresponding to the element $w_{{m,i}%
}\in W(P;G)$. With respect to the Schubert basis $\{s_{{m,i}}\mid 1\leq
i\leq $\textsl{\ }$\beta (m)\}$ on $H^{2m}(G/P)$ every monomial $y^{\alpha
}\in B(2m)$ has the unique expansion

\begin{quote}
$\pi _{{m}}(y^{\alpha })=c_{{\alpha ,1}}\cdot s_{{m,1}}+\cdots +c_{{\alpha
,\beta (m)}}\cdot s_{{m,\beta (m)}}$, $c_{{\alpha ,i}}\in \mathbb{Z}$,
\end{quote}

\noindent where the coefficients $c_{{\alpha ,i}}$ can be evaluated by the 
\textsl{Characteristics}. The matrix $M(f_{{m}})=(c_{{\alpha ,i}})_{{%
b(2m)\times \beta (m)}}$ so obtained is called \textsl{the structure matrix
of }$f_{{m}}$. The built-in function \textsl{Null-space} in \textsl{%
Mathematica} transforms $M(f_{{m}})$ to another matrix $N(f_{{m}})$ in the
fashion

\begin{quote}
\textsl{In:}=Null-space[$M(f_{{m}})$]

\textsl{Out:}= a matrix\textsl{\ }$N(f_{{m}})=(b_{{j,\alpha }})_{{%
(b(2m)-\beta (m))\times b(2m)}}$,
\end{quote}

\noindent whose significance is shown by the following fact.

\bigskip

\noindent \textbf{Lemma 4.7.} \textsl{The set }$\kappa _{{i}%
}=\sum_{y^{\alpha }\in B(2m)}b_{{i,\alpha }}\cdot y^{\alpha }$, $1\leq i\leq
b(2m)-\beta (m)$, \textsl{of polynomials is a basis of the }$\mathbb{Z}$ 
\textsl{module} $\ker f_{{m}}$.\hfill $\square$

\bigskip

\textbf{4.4.3. The Giambelli polynomials (i.e. the Schubert polynomials \cite%
{Bi}) }For the unitary group $G=U(n)$ of rank $n$ with parabolic subgroup $%
P=U(k)\times U(n-k)$ the flag manifold $G_{{n,k}}=G/P$ is the \textsl{%
Grassmannian} of $k$-planes through the origin on $\mathbb{C}^{n}$. Let $%
1+c_{1}+\cdots +c_{k}$ be the total Chern class of the canonical $k$-bundle
on $G_{{n,k}}$. Then the $c_{i}$'s can be identified with appropriate
Schubert classes on $G_{{n,k}}$ (i.e. \textsl{the} \textsl{special Schubert
class} on $G_{{n,k}}$), and one has \textsl{the classical Schubert
presentation}

\begin{quote}
$H^{\ast }(G_{{n,k}})=\mathbb{Z}[c_{{1}},\ldots ,c_{{k}}]/\left\langle r_{{%
n-k+1}},\ldots ,r_{{n}}\right\rangle $,
\end{quote}

\noindent where $r_{{j}}$ is the component of the formal inverse of $1+c_{{1}%
}+\cdots +c_{{k}}$ in degree $j$. It follows that every Schubert class $%
s_{w} $ on $G_{{n,k}}$ can be written as a polynomial $\mathcal{G}_{{w}}(c_{{%
1}},\ldots ,c_{{k}})$ in the special ones, and such an expression is
afforded by the classical \textsl{Giambelli formula} \cite[p.112]{Hil}.

In general, assume that $G/P$ is a flag variety, and that a Schubert
presentation (4.1) of the ring $H^{\ast }(G/P)$ has been specified. Then
each Schubert class $s_{{w}}$ of $G/P$ can be expressed as a polynomial $%
\mathcal{G}_{{w}}(y_{{1}},\ldots ,y_{{k}})$ in these special ones, and such
an expression will be called a \textsl{Giambelli polynomial} of the class $%
s_{{w}}$. Based on the \textsl{Characteristics} a program implementing $%
\mathcal{G}_{{w}}(y_{{1}},\ldots ,y_{{k}})$ has been compiled, whose
function is summarized below.

\begin{quote}
\textbf{Algorithm 4.8:} \textsl{Giambelli polynomials}

\textbf{Input:} \textsl{A set }$\{y_{{1}},\ldots ,y_{{k}}\}$\textsl{\ of
special Schubert classes on }$G/P$\textsl{.}

\textbf{Output:} \textsl{Giambelli polynomials }$\mathcal{G}_{{w}}(y_{{1}%
},\ldots ,y_{{k}})$\textsl{\ for all }$w\in W(P;G)$.
\end{quote}

We clarify the details in this program. By (4.13) we can write the ordered
monomial basis $B(2m)$ of\textsl{\ }$\mathbb{Z}[y_{{1}},\ldots ,y_{{k}%
}]^{(2m)}$ as $\{y^{\alpha _{{1}}},\ldots ,y^{\alpha _{{b(2m)}}}\}$. The
corresponding structure matrix $M(f_{{m}})$ in degree $2m$ then satisfies

\begin{quote}
$\left( 
\begin{array}{c}
y^{\alpha _{{1}}} \\ 
\vdots \\ 
y^{\alpha _{{b(2m)}}}%
\end{array}%
\right) =M(f_{{m}})\left( 
\begin{array}{c}
s_{{m,1}} \\ 
\vdots \\ 
s_{{m,\beta (m)}}%
\end{array}%
\right) $.
\end{quote}

\noindent Since $f_{{m}}$ surjects the matrix $M(f_{{m}})$ has a $\beta
(m)\times \beta (m)$ minor equal to $\pm 1$. The standard integral row and
column operation diagonalizing $M(f_{{m}})$ \cite[p.162-164]{Sch} provides
us with two unique invertible matrices $P=P_{{b(2m)\times b(2m)}}$ and $Q=Q_{%
{\beta (m)\times \beta (m)}}$ that satisfy

\begin{enumerate}
\item[(4.14)] $P\cdot M(f_{{m}})\cdot Q=\left( 
\begin{array}{c}
I_{{\beta (m)}} \\ 
C%
\end{array}%
\right) _{{b(2m)\times \beta (m)}}$,
\end{enumerate}

\noindent where $I_{{\beta (m)}}$ is the identity matrix of rank $\beta (m)$%
. The \textsl{Giambelli polynomials} is realized by the procedure below.

\begin{quote}
\textbf{Step 1.} Compute $M(f_{{m}})$ using the \textsl{Characteristics};

\textbf{Step 2.} Diagonalize $M(f_{{m}})$ to get the matrices $P$ and $Q$;

\textbf{Step 3.} Set $\left( 
\begin{array}{c}
\mathcal{G}_{{m,1}} \\ 
\vdots \\ 
\mathcal{G}_{{m,\beta (m)}}%
\end{array}%
\right) =Q\cdot \lbrack P]\left( 
\begin{array}{c}
y^{\alpha _{{1}}} \\ 
\vdots \\ 
y^{\alpha _{{b(2m)}}}%
\end{array}%
\right) $,
\end{quote}

\noindent where $[P]$ is formed by the first $\beta (m)$ rows of $P$.
Obviously, the polynomial $\mathcal{G}_{{m,j}}$ so obtained depends only on
the special Schubert classes $\{y_{{1}},\ldots ,$ $y_{{k}}\}$ on $G/P$, and
is a Giambelli polynomial of\textsl{\ }$s_{m,j}$, $1\leq j\leq {\beta (m)}$.

\section{Application to the flag manifolds $G/T$}

\begin{quote}
\textsl{A calculus, or science of calculation, is one which has organized
processes by which passage is made, mechanically, from one result to
another. }-De Morgan.
\end{quote}

Among all the flag manifolds $G/P$ associated to a Lie group $G$ it is the
complete flag manifold $G/T$ that is of crucial importance, since the
inclusion $T\subseteq P\subset G$ of subgroups induces the fibration

\begin{quote}
$P/T\hookrightarrow G/T\overset{\pi }{\longrightarrow }G/P$
\end{quote}

\noindent in which the induced map $\pi ^{\ast }$ embeds the ring $H^{\ast
}(G/P)$ as a subring of $H^{\ast }(G/T)$, see the proof of Theorem 4.3.
Further, according to E. Cartan \cite[p.674]{Wh} all the $1$-connected
simple Lie groups consist of the three infinite families $SU(n)$, $Sp(n)$, $%
Spin(n)$ of the classical groups, as well as the five exceptional ones: $%
G_{2},F_{4},E_{6},E_{7},E_{8}$, while for any compact connected Lie group $G$
with a maximal torus $T$ one has a diffeomorphism

\begin{quote}
$G/T=G_{1}/T_{1}\times \cdots \times G_{k}/T_{k}$ \noindent
\end{quote}

\noindent with each $G_{i}$ an $1$-connected simple Lie group and $%
T_{i}\subset G_{i}$ a maximal torus. Thus, the problem of finding Schubert
presentations of flag manifolds may be reduced to the special cases $G/T$
where $G$ is $1$-connected and simple.

In this section we determine the Schubert presentation of the ring $H^{\ast
}(G/T)$ in accordance to $G$ is classical or exceptional. Recall that for a
simple Lie group $G$ with rank $n$ the fundamental dominant weights $\left\{
\omega _{1},\ldots ,\omega _{n}\right\} $ of $G$ \cite{BH} is precisely the
Schubert basis of $H^{2}(G/T)$ \cite[Lemma 2.4]{Du}.

\subsection{The ring $H^{\ast }(G/T)$ for a classical $G$}

If $G=SU(n+1)$ or $Sp(n)$ Borel \cite{B1} has shown that

\begin{quote}
$H^{\ast }(G/T)=\mathbb{Z}\left[ \omega _{1},\ldots ,\omega _{n}\right]
/\left\langle \mathbb{Z}\left[ \omega _{1},\ldots ,\omega _{n}\right]
^{+,W}\right\rangle $,
\end{quote}

\noindent where $\mathbb{Z}\left[ \omega _{1},\cdots ,\omega _{n}\right]
^{+,W}$\ is the ring of the integral Weyl invariants of $G$ in positive
degrees. It follows that if we let $c_{r}(G)\in H^{2r}(G/T)$ be respectively
the $r^{th}$ element symmetric polynomial in the sets

\begin{quote}
$\left\{ \omega _{{\footnotesize 1}},\omega _{{\footnotesize k}}-\omega _{%
{\footnotesize k-1}},-\omega _{{\footnotesize n}}\mid \text{ }2\leq k\leq
n\right\} $ or $\left\{ \pm \omega _{{\footnotesize 1}},\pm (\omega _{%
{\footnotesize k}}-\omega _{{\footnotesize k-1}})\mid 2\leq k\leq n\right\} $%
,
\end{quote}

\noindent then we have

\bigskip

\noindent \textbf{Theorem 5.1.} \textsl{For }$G=SU(n)$\textsl{\ or }$Sp(n)$ 
\textsl{Schubert presentation of }$H^{\ast }(G/T)$\textsl{\ is}

\begin{enumerate}
\item[(5.1)] $H^{\ast }(SU(n)/T)=\mathbb{Z}\left[ \omega _{1},\ldots ,\omega
_{n-1}\right] /\left\langle c_{2},\ldots ,c_{n}\right\rangle $\textsl{,} $%
c_{r}=c_{r}(SU(n))$\textsl{;}

\item[(5.2)] $H^{\ast }(Sp(n)/T)=\mathbb{Z}\left[ \omega _{1},\ldots ,\omega
_{n}\right] /\left\langle c_{2},\ldots ,c_{2n}\right\rangle $\textsl{,} $%
c_{2r}=c_{2r}(Sp(n))$\textsl{.}\hfill $\square $
\end{enumerate}

Turning to the group $G=Spin(2n)$ let $y_{k}$ be the Schubert class on $%
Spin(2n)/T$ associated to the Weyl group element

\begin{quote}
$w_{k}=\sigma \lbrack n-k,\ldots ,n-2,n-1]$, $2\leq k\leq n-1$
\end{quote}

\noindent (in the notation of Section 4.1). According to Marlin \cite[%
Proposition 3]{M}

\begin{enumerate}
\item[(5.3)] $H^{\ast }(Spin(2n)/T)=\mathbb{Z}[\omega _{1},\cdots ,\omega
_{n},y_{2},\cdots ,y_{n-1}]/\left\langle \delta _{i},\xi _{j},\mu
_{k}\right\rangle $
\end{enumerate}

\noindent where

\begin{quote}
$\delta _{i}:=2y_{i}-c_{i}(\omega _{1},\ldots ,\omega _{n})$, $1\leq i\leq
n-1$,

$\xi _{j}:=y_{2j}+(-1)^{j}y_{j}^{2}+2\sum\limits_{1\leq r\leq
j-1}(-1)^{r}y_{r}y_{2j-r}$, $1\leq j\leq \left[ \frac{n-1}{2}\right] $,

$\mu _{k}:=(-1)^{k}y_{k}^{2}+2\sum\limits_{2k-n+1\leq r\leq
k-1}(-1)^{r}y_{r}y_{2k-r}$, $\left[ \frac{n}{2}\right] \leq k\leq n-1$,
\end{quote}

\noindent and where $c_{i}(\omega _{1},\ldots ,\omega _{n})$ is the $i^{th}$
elementary symmetric function on set

\begin{quote}
$\left\{ \omega _{n},\omega _{i}-\omega _{i-1},\omega _{n-1}+\omega
_{n}-\omega _{n-2},\omega _{n-1}-\omega _{n}\mid 2\leq i\leq n-2\right\} $.
\end{quote}

\noindent Since each Schubert class $y_{2j}$ with $1\leq j\leq \left[ \frac{%
n-1}{2}\right] $ can be expressed as a polynomial in the $y_{2i+1}$'s by the
relations of the type $\xi _{k}$, we obtain that

\bigskip

\noindent \textbf{Theorem 5.2.} \textsl{The Schubert presentation of }$%
H^{\ast }(Spin(2n)/T)$\textsl{\ is}

\begin{enumerate}
\item[(5.4)] $H^{\ast }(Spin(2n)/T)=\mathbb{Z}[\omega _{1},\ldots ,\omega
_{n},y_{3},y_{5},\ldots ,y_{2\left[ \frac{n-1}{2}\right] -1}]/\left\langle
r_{i},h_{k}\right\rangle $,
\end{enumerate}

\noindent \textsl{where }$r_{i}$\textsl{\ and }$h_{k}$\textsl{\ are the
polynomials obtained respectively from }$\delta _{i}$\textsl{\ and }$\mu
_{k} $\textsl{\ by replacing the classes }$y_{2r}$\textsl{\ with the
polynomials (by the relation }$\xi _{r}$\textsl{)}

\begin{quote}
$(-1)^{r-1}y_{r}^{2}+2\sum\limits_{1\leq k\leq r-1}(-1)^{k-1}y_{k}y_{2r-k}$%
\textsl{.}\hfill $\square$
\end{quote}

Similarly, if $G=Spin(2n+1)$ one can deduce a Schubert presentation of the
ring $H^{\ast }(G/T)$ from Marlin's formula \cite[Proposition 2]{M}.

\bigskip

\noindent \textbf{Remark 5.3.} For the classical Lie groups $G$ the
Giambelli polynomials (i.e. \textsl{Schubert polynomials}) of the Schubert
classes on $G/T$ have been determined by Billey and Haiman \cite{Bi}.

In comparison with Marlin's formula (5.3) the presentation (5.4) is more
concise for involving fewer generators and relations. A basic requirement of
topology is to present the cohomology of a space $X$ by a minimal system of
generators and relations, so that the (rational) minimal model and $\kappa $%
-invariants of the Postnikov tower of $X$ \cite{GM} can be formulated
accordingly.\hfill $\square $

\subsection{The ring $H^{\ast }(G/T)$ for an exceptional $G$}

Having clarified the Schubert presentation of the ring $H^{\ast }(G/T)$ for
the classical $G$ we proceed to the exceptional cases $G=F_{4},E_{6}$ or $%
E_{7}$ (the result for the case $E_{8}$ comes from the same calculation,
only the presentation \cite[Theorem 5.1]{DZ3} is slightly lengthy). In these
cases the dimension $s=\dim G/T$ and the number $t$ of the Schubert classes
on $G/T$ are

\begin{quote}
$(s,t)=(48,1152),(72,51840)$ or $(126,2903040)$,
\end{quote}

\noindent respectively. Instead of describing the ring $H^{\ast }(G/T)$
using the totality of $t^{3}$ Littlewood-Richardson coefficients $%
c_{u,v}^{w} $ (with $c_{u,v}^{w}=0$ for $l(w)\neq l(u)+l(v)$ being
understood) the idea of Schubert presentation brings us the following
concise and explicit formulae of the ring $H^{\ast }(G/T)$.

\bigskip

\noindent \textbf{Theorem 5.4.} \textsl{For }$G=F_{4},E_{6}$\textsl{\ and }$%
E_{7}$\textsl{\ the Schubert presentations of the cohomologies }$H^{\ast
}(G/T)$\textsl{\ are}

\begin{enumerate}
\item[(5.5)] $H^{\ast }(F_{4}/T)=\mathbb{Z}[\omega _{1},\omega _{2},\omega
_{3},\omega _{4},y_{3},y_{4}]/\left\langle \rho _{2},\rho
_{4},r_{3},r_{6},r_{8},r_{12}\right\rangle $\textsl{, where}

$\rho _{2}=c_{2}-4\omega _{1}^{2}$;

$\rho _{4}=3y_{4}+2\omega _{1}y_{3}-c_{4}$;

$r_{3}=2y_{3}-\omega _{1}^{3}$;

$r_{6}=y_{3}^{2}+2c_{6}-3\omega _{1}^{2}y_{4}$;

$r_{8}=3y_{4}^{2}-\omega _{1}^{2}c_{6}$;

$r_{12}=y_{4}^{3}-c_{6}^{2}$.

\item[(5.6)] $H^{\ast }(E_{6}/T)=\mathbb{Z}[\omega _{1},\ldots ,\omega
_{6},y_{3},y_{4}]/\left\langle \rho _{2},\rho _{3},\rho _{4},\rho
_{5},r_{6},r_{8},r_{9},r_{12}\right\rangle $\textsl{, where}

$\rho _{2}=4\omega _{2}^{2}-c_{2}$;

$\rho _{3}=2y_{3}+2\omega _{2}^{3}-c_{3}$;

$\rho _{4}=3y_{4}+\omega _{2}^{4}-c_{4}$;

$\rho _{5}=2\omega _{2}^{2}y_{3}-\omega _{2}c_{4}+c_{5}$;

$r_{6}=y_{3}^{2}-\omega _{2}c_{5}+2c_{6}$;

$r_{8}=3y_{4}^{2}-2c_{5}y_{3}-\omega _{2}^{2}c_{6}+\omega _{2}^{3}c_{5};$

$r_{9}=2y_{3}c_{6}-\omega _{2}^{3}c_{6}$;

$r_{12}=y_{4}^{3}-c_{6}^{2}$.

\item[(5.7)] $H^{\ast }(E_{7}/T)=\mathbb{Z}[\omega _{1},\ldots ,\omega
_{7},y_{3},y_{4},y_{5},y_{9}]/\left\langle \rho _{i},r_{j}\right\rangle $%
\textsl{, where}

$\rho _{2}=4\omega _{2}^{2}-c_{2}$;

$\rho _{3}=2y_{3}+2\omega _{2}^{3}-c_{3}$;

$\rho _{4}=3y_{4}+\omega _{2}^{4}-c_{4}$;

$\rho _{5}=2y_{5}-2\omega _{2}^{2}y_{3}+\omega _{2}c_{4}-c_{5}$;

$r_{6}=y_{3}^{2}-\omega _{2}c_{5}+2c_{6}$;

$r_{8}=3y_{4}^{2}+2y_{3}y_{5}-2y_{3}c_{5}+2\omega _{2}c_{7}-\omega
_{2}^{2}c_{6}+\omega _{2}^{3}c_{5}$;

$r_{9}=2{y_{9}}+2{y_{4}y_{5}}-2{y_{3}c_{6}}-{\omega _{2}^{2}c_{7}}+{\omega
_{2}^{3}c_{6}}$;

$r_{10}=y_{5}^{2}-2y_{3}c_{7}+\omega _{2}^{3}c_{7}$;

$r_{12}=y_{4}^{3}-4y_{5}c_{7}-c_{6}^{2}-2y_{3}y_{9}-2y_{3}y_{4}y_{5}+2\omega
_{2}y_{5}c_{6}+3\omega _{2}y_{4}c_{7}+c_{5}c_{7}$;

$r_{14}=c_{7}^{2}-2y_{5}y_{9}+2y_{3}y_{4}c_{7}-\omega _{2}^{3}y_{4}c_{7}$;

$%
r_{18}=y_{9}^{2}+2y_{5}c_{6}c_{7}-y_{4}c_{7}^{2}-2y_{4}y_{5}y_{9}+2y_{3}y_{5}^{3}-5\omega _{2}y_{5}^{2}c_{7} 
$,
\end{enumerate}

\noindent \textsl{where the }$c_{r}$\textsl{'s are the polynomials }$%
c_{r}(P) $ \textsl{in the weights }$\omega _{1},\cdots ,\omega _{n}$ \textsl{%
defined in (5.17); and where the }$y_{i}$\textsl{'s are the Schubert classes
on }$G/T $\textsl{\ associated to the Weyl group elements tabulated below:}

\begin{enumerate}
\item[(5.8)] 
\begin{tabular}{|l|l|l|l|}
\hline
$y_{i}$ & $F_{4}/T$ $(F/P_{\{1\}})$ & $E_{6}/T$ $(E_{6}/P_{\{2\}})$ & $%
E_{7}/T$ $(E_{7}/P_{\{2\}})$ \\ \hline\hline
$y_{3}$ & $\sigma _{\lbrack 3,2,1]}$ & $\sigma _{\lbrack 5,4,2]}$ & $\sigma
_{\lbrack 5,4,2]}$ \\ 
$y_{4}$ & $\sigma _{\lbrack 4,3,2,1]}$ & $\sigma _{\lbrack 6,5,4,2]}$ & $%
\sigma _{\lbrack 6,5,4,2]}$ \\ 
$y_{5}$ &  &  & $\sigma _{\lbrack {7,6,5,4,2}]}$ \\ 
$y_{6}$ & $\sigma _{\lbrack 3,2,4,3,2,1]}$ & $\sigma _{\lbrack 1,3,6,5,4,2]}$
& $\sigma _{\lbrack 1,3,6,5,4,2]}$ \\ 
$y_{7}$ &  &  & $\sigma _{\lbrack 1,3,7,6,5,4,2]}$ \\ 
$y_{9}$ &  &  & $\sigma _{\lbrack 1,5,4,3,7,6,5,4,2]}$. \\ \hline
\end{tabular}
\end{enumerate}

\bigskip

To reduce the computational complexity of deriving Theorem 5.4 we choose for
each $G=F_{4},E_{6}$ or $E_{7}$ a parabolic subgroup $P$ associated to a
singleton $K=\{i\}$, where the index $i$, as well as the isomorphism types
of $P$ and its simple part $P^{s}$, is stated in the table below

\begin{enumerate}
\item[(5.9)] 
\begin{tabular}{l|l|l|l}
\hline
$G$ & $F_{4}$ & $E_{6}$ & $E_{7}$ \\ \hline\hline
$i$ & $1$ & $2$ & $2$ \\ \hline
$P;P^{s}$ & $Sp(3)\cdot S^{1};Sp(3)$ & $SU(6)\cdot S^{1};SU(6)$ & $%
SU(7)\cdot S^{1};SU(7)$ \\ \hline
\end{tabular}%
.
\end{enumerate}

\noindent In view of the circle bundle associated to $P$ (see in Section 4.3)

\begin{enumerate}
\item[(5.10)] $S^{1}\hookrightarrow G/P^{s}\overset{p}{\rightarrow }G/P$
\end{enumerate}

\noindent the calculation will be divided into three steps, in accordance to
cohomologies of the three homogeneous spaces $G/P^{s}$, $G/P$ and $G/T$.

\bigskip

\noindent \textbf{Step 1. The cohomologies }$H^{\ast }(G/P^{s})$\textbf{.}
By the formulae (4.8) and (4.9) the additive structure of $H^{\ast
}(G/P^{s}) $ is determined by the homomorphisms

\begin{quote}
$H^{2r-2}(G/P)\overset{\cup \omega }{\rightarrow }H^{2r}(G/P)$.
\end{quote}

\noindent Explicitly, with respect to the Schubert basis $\{s_{{r,1}},\ldots
,$ $s_{{r,\beta (r)}}\}$ of $H^{2r}(G/P)$, $\beta (r)=\left\vert
W^{r}(P;G)\right\vert $, the expansions

\begin{quote}
$\omega \cup s_{{r-1,i}}=\sum a_{{i,j}}\cdot s_{{r,j}}$
\end{quote}

\noindent give rise to a ${\beta (r-1)\times \beta (r)}$ matrix $A_{r}$ that
satisfies the linear system

\begin{quote}
$%
\begin{array}{rcl}
\left( 
\begin{array}{l}
\omega \cup s_{{r-1,1}} \\ 
\omega \cup s_{{r-1,2}} \\ 
\text{ \ \ \ }\vdots \\ 
\omega \cup s_{{r-1,\beta (r-1)}}%
\end{array}%
\right) & = & A_{{r}}\left( 
\begin{array}{l}
s_{{r,1}} \\ 
s_{{r,2}} \\ 
\text{ }\vdots \\ 
s_{{r,\beta (r)}}%
\end{array}%
\right)%
\end{array}%
$.
\end{quote}

\noindent Since $\omega \cup s_{{r-1,i}}$ is a monomial in Schubert classes,
the \textsl{Characteristics }is applicable to evaluate\textsl{\ }the entries
of $A_{{r}}$. Diagonalizing $A_{{r}}$ by the integral row and column
reductions (\cite[p.162-166]{S}) one obtains the non-trivial groups $%
H^{r}(G/P^{s})$, together with their basis, as that tabulated below, where

\begin{quote}
1) $\overline{y}_{i}:=p^{\ast }(y_{i})$ with $y_{i}$ the Schubert classes in
table (5.8);

2) For simplicity the non-trivial groups $H^{r}(G/P^{s})$ are printed only
up to the stage where all the generators and relations of the ring $%
H^{even}(G/P^{s})$ emerge.
\end{quote}

\noindent \textbf{Table 1. Non-trivial cohomologies of }$F_{4}/P^{s}:$

\begin{center}
\begin{tabular}{|l|l|l|}
\hline
nontrivial $H^{k}$ & basis elements & relations \\ \hline\hline
$H^{6}\cong \mathbb{Z}_{{2}}$ & $\bar{s}_{{3,1}}(=\overline{y}_{3})$ & $2%
\overline{y}_{3}=0$ \\ \hline
$H^{8}\cong \mathbb{Z}$ & $\bar{s}_{{4,2}}(=\overline{y}_{4})$ &  \\ \hline
$H^{12}\cong \mathbb{Z}_{{4}}$ & $\bar{s}_{{6,2}}(=\overline{y}_{6})$ & $-2%
\overline{y}_{6}=\overline{y}_{3}^{2}$ \\ \hline
$H^{14}\cong \mathbb{Z}_{{2}}$ & $\overline{y}_{3}\overline{y}_{4}$ &  \\ 
\hline
$H^{16}\cong \mathbb{Z}_{{3}}$ & $\overline{y}_{4}^{2}$ & $3\overline{y}%
_{4}^{2}=0$ \\ \hline
$H^{18}\cong \mathbb{Z}_{{2}}$ & $\overline{y}_{3}\overline{y}_{6}$ &  \\ 
\hline
$H^{20}\cong \mathbb{Z}_{{4}}$ & $\overline{y}_{4}\overline{y}_{6}$ &  \\ 
\hline
$H^{26}\cong \mathbb{Z}_{{2}}$ & $\overline{y}_{3}\overline{y}_{4}\overline{y%
}_{6}$ &  \\ \hline
$H^{23}\cong \mathbb{Z}$ & $d_{{23}}=\beta ^{-1}(2s_{{11,1}}-s_{{11,2}})$ \
\ \  &  \\ \hline
\end{tabular}
\end{center}

\noindent \textbf{Table 2. Non-trivial cohomologies of }$E_{6}/P^{s}:$

\begin{center}
\begin{tabular}{|l|l|l|}
\hline
nontrivial $H^{k}$ & basis elements & relations \\ \hline\hline
$H^{6}\cong \mathbb{Z}$ & $\bar{s}_{{3,2}}(=\overline{y}_{3})$ &  \\ \hline
$H^{8}\cong \mathbb{Z}$ & $\bar{s}_{{4,3}}(=\overline{y}_{4})$ &  \\ \hline
$H^{12}\cong \mathbb{Z}$ & $\bar{s}_{{6,1}}(=\overline{y}_{6})$ & $-2%
\overline{y}_{6}=\overline{y}_{3}^{2}$ \\ \hline
$H^{14}\cong \mathbb{Z}$ & $\overline{y}_{3}\overline{y}_{4}$ &  \\ \hline
$H^{16}\cong \mathbb{Z}_{{3}}$ & $\overline{y}_{4}^{2}$ & $3\overline{y}%
_{4}^{2}=0$ \\ \hline
$H^{18}\cong \mathbb{Z}_{{2}}$ & $\overline{y}_{3}\overline{y}_{6}$ & $2%
\overline{y}_{3}\overline{y}_{6}=0$ \\ \hline
$H^{20}\cong \mathbb{Z}$ & $\overline{y}_{4}\overline{y}_{6}$ &  \\ \hline
$H^{22}\cong \mathbb{Z}_{{3}}$ & $\overline{y}_{3}\overline{y}_{4}^{2}$ & 
\\ \hline
$H^{26}\cong \mathbb{Z}_{{2}}$ & $\overline{y}_{3}\overline{y}_{4}\overline{y%
}_{6}$ &  \\ \hline
$H^{28}\cong \mathbb{Z}_{{3}}$ & $\overline{y}_{4}^{2}\overline{y}_{6}$ & 
\\ \hline
$H^{23}\cong \mathbb{Z}$ & ${\small d}_{{\small 23}}{\small =\beta }^{%
{\small -1}}{\small (s}_{{\small 11,1}}{\small -s}_{{\small 11,2}}{\small -s}%
_{{\small 11,3}}{\small +s}_{{\small 11,4}}{\small -s}_{{\small 11,5}}%
{\small +s}_{{\small 11,6}}{\small )}$ &  \\ \hline
$H^{29}\cong \mathbb{Z}$ & ${\small d}_{{\small 29}}{\small =\beta }^{%
{\small -1}}{\small (-s}_{{\small 14,1}}{\small +s}_{{\small 14,2}}{\small +s%
}_{{\small 14,4}}{\small -\ s}_{{\small 14,5}}{\small )}$ & $2d_{{29}}=\pm 
\overline{y}_{3}d_{{23}}$ \\ \hline
\end{tabular}
\end{center}

\noindent \textbf{Table 3. Non-trivial cohomologies of }$E_{7}/P^{s}$

\begin{center}
\begin{tabular}{|l|l|}
\hline
nontrivial $H^{k}$ & basis elements \\ \hline\hline
$H^{6}\cong \mathbb{Z}$ & $s_{3,2}=\bar{y}_{3}$ \\ \hline
$H^{8}\cong \mathbb{Z}$ & $s_{4,3}=\bar{y}_{4}$ \\ \hline
$H^{10}\cong \mathbb{Z}$ & $s_{5,4}=\bar{y}_{5}$ \\ \hline
$H^{12}\cong \mathbb{Z}$ & $s_{6,5}=\bar{y}_{3}^{2}+\bar{y}_{6}$ \\ \hline
$H^{14}\cong \mathbb{Z}\oplus \mathbb{Z}$ & $s_{7,6}=-\bar{y}_{7};$ $\bar{y}%
_{3}\bar{y}_{4}$ \\ \hline
$H^{16}\cong \mathbb{Z}$ & $\bar{y}_{3}\bar{y}_{5}-2\bar{y}_{4}^{2}$ \\ 
\hline
$H^{18}\cong \mathbb{Z}_{2}\oplus \mathbb{Z}\oplus \mathbb{Z}$ & $\bar{y}%
_{3}^{3}+\bar{y}_{3}\bar{y}_{6}+\bar{y}_{4}\bar{y}_{5}+\bar{y}_{9};-\bar{y}%
_{4}\bar{y}_{5};\bar{y}_{3}^{3}+\bar{y}_{9}$ \\ \hline
$H^{20}\cong \mathbb{Z}\oplus \mathbb{Z}$ & $\bar{y}_{3}\bar{y}_{7}-\bar{y}%
_{5}^{2};$ $\bar{y}_{3}\bar{y}_{7}+\bar{y}_{4}\bar{y}_{6}$ \\ \hline
$H^{22}\cong \mathbb{Z}\oplus \mathbb{Z}$ & $\bar{y}_{3}^{2}\bar{y}_{5}-\bar{%
y}_{3}\bar{y}_{4}^{2}+\bar{y}_{4}\bar{y}_{7};$ $-2\bar{y}_{4}\bar{y}_{7}+%
\bar{y}_{5}\bar{y}_{6}$ \\ \hline
$H^{24}\cong \mathbb{Z}_{2}\oplus \mathbb{Z}\oplus \mathbb{Z}$ & $\bar{y}%
_{3}^{4}+\bar{y}_{3}^{2}\bar{y}_{6}+\bar{y}_{3}\bar{y}_{4}\bar{y}_{5}+\bar{y}%
_{3}\bar{y}_{9}$ \\ 
& $\bar{y}_{3}^{2}\bar{y}_{6}+\bar{y}_{3}\bar{y}_{4}\bar{y}_{5}-\bar{y}%
_{4}^{3};$ $-\bar{y}_{4}^{3}+\bar{y}_{5}\bar{y}_{7}+\bar{y}_{6}^{2}$ \\ 
\hline
$H^{26}\cong \mathbb{Z}_{2}\oplus \mathbb{Z}\oplus \mathbb{Z}$ & $\bar{y}%
_{3}^{3}\bar{y}_{4}+\bar{y}_{3}\bar{y}_{4}\bar{y}_{6}+\bar{y}_{4}^{2}\bar{y}%
_{5}+\bar{y}_{4}\bar{y}_{9}$ \\ 
& $\bar{y}_{3}\bar{y}_{4}\bar{y}_{6}+\bar{y}_{3}\bar{y}_{5}^{2}-3\bar{y}%
_{4}^{2}\bar{y}_{5}+\bar{y}_{6}\bar{y}_{7}$ \\ 
& $\bar{y}_{3}^{3}\bar{y}_{4}+\bar{y}_{3}^{2}\bar{y}_{7}+2\bar{y}_{4}^{2}%
\bar{y}_{5}+\bar{y}_{4}\bar{y}_{9}$ \\ \hline
$H^{28}\cong \mathbb{Z}_{2}\oplus \mathbb{Z}\oplus \mathbb{Z}$ & $\bar{y}%
_{3}^{3}\bar{y}_{5}+\bar{y}_{3}\bar{y}_{5}\bar{y}_{6}+\bar{y}_{4}\bar{y}%
_{5}^{2}+\bar{y}_{5}\bar{y}_{9}$ \\ 
& $3\bar{y}_{3}^{2}\bar{y}_{4}^{2}+5\bar{y}_{4}\bar{y}_{5}^{2}+\bar{y}_{5}%
\bar{y}_{9}$ \\ 
& $\bar{y}_{3}^{3}\bar{y}_{5}+2\bar{y}_{3}^{2}\bar{y}_{4}^{2}+\bar{y}_{3}%
\bar{y}_{4}\bar{y}_{7}+\bar{y}_{3}\bar{y}_{5}\bar{y}_{6}+4\bar{y}_{4}\bar{y}%
_{5}^{2}+\bar{y}_{5}\bar{y}_{9}$ \\ \hline
$H^{30}\cong \mathbb{Z}_{2}\oplus \mathbb{Z}\oplus \mathbb{Z}$ & $-\bar{y}%
_{3}^{5}-\bar{y}_{3}\bar{y}_{4}^{3}+\bar{y}_{4}\bar{y}_{5}\bar{y}_{6}+\bar{y}%
_{5}^{3}+\bar{y}_{6}\bar{y}_{9}$ \\ 
& $-\bar{y}_{3}^{5}+\bar{y}_{3}^{2}\bar{y}_{4}\bar{y}_{5}-3\bar{y}_{3}\bar{y}%
_{4}^{3}+\bar{y}_{3}\bar{y}_{5}\bar{y}_{7}-\bar{y}_{4}^{2}\bar{y}_{7}+\bar{y}%
_{5}^{3}$ \\ 
& $\bar{y}_{3}^{2}\bar{y}_{4}\bar{y}_{5}+\bar{y}_{3}\bar{y}_{5}\bar{y}_{7}+%
\bar{y}_{4}^{2}\bar{y}_{7}+\bar{y}_{4}\bar{y}_{5}\bar{y}_{6}$ \\ \hline
$H^{32}\cong \mathbb{Z}_{2}\oplus \mathbb{Z}_{2}\oplus \mathbb{Z}\oplus 
\mathbb{Z}$ & $\bar{y}_{3}^{4}\bar{y}_{4}+\bar{y}_{3}^{2}\bar{y}_{4}\bar{y}%
_{6}+\bar{y}_{3}\bar{y}_{4}^{2}\bar{y}_{5}+\bar{y}_{3}\bar{y}_{4}\bar{y}_{9}$
\\ 
& $\bar{y}_{3}^{2}\bar{y}_{5}^{2}+\bar{y}_{3}\bar{y}_{6}\bar{y}_{7}+\bar{y}%
_{4}\bar{y}_{5}\bar{y}_{7}+\bar{y}_{5}^{2}\bar{y}_{6}+\bar{y}_{7}\bar{y}_{9}$
\\ 
& $\bar{y}_{3}^{4}\bar{y}_{4}+\bar{y}_{3}\bar{y}_{4}^{2}\bar{y}_{5}+\bar{y}%
_{3}\bar{y}_{4}\bar{y}_{9}+\bar{y}_{3}\bar{y}_{6}\bar{y}_{7}+\bar{y}_{4}^{4}$
\\ 
& $\bar{y}_{3}^{2}\bar{y}_{4}\bar{y}_{6}+\bar{y}_{3}\bar{y}_{4}^{2}\bar{y}%
_{5}-\bar{y}_{4}^{4}+\bar{y}_{4}\bar{y}_{5}\bar{y}_{7}+\bar{y}_{5}^{2}\bar{y}%
_{6}$ \\ \hline
$H^{34}\cong \mathbb{Z}_{38}\oplus \mathbb{Z}$ & $\bar{y}_{3}^{4}\bar{y}_{5}+%
\bar{y}_{3}\bar{y}_{4}\bar{y}_{5}^{2}+\bar{y}_{3}\bar{y}_{5}\bar{y}_{9}+2%
\bar{y}_{4}^{3}\bar{y}_{5}$ \\ 
& $\bar{y}_{3}^{4}\bar{y}_{5}+2\bar{y}_{3}^{3}\bar{y}_{4}^{2}+\bar{y}_{3}^{2}%
\bar{y}_{4}\bar{y}_{7}+9\bar{y}_{4}^{3}\bar{y}_{5}+\bar{y}_{4}\bar{y}_{6}%
\bar{y}_{7}$ \\ \hline
$H^{36}\cong \mathbb{Z}_{2}\oplus \mathbb{Z}_{2}\oplus \mathbb{Z}_{20}\oplus 
\mathbb{Z}$ & $\bar{y}_{3}^{3}\bar{y}_{4}\bar{y}_{5}+\bar{y}_{3}\bar{y}_{4}%
\bar{y}_{5}\bar{y}_{6}+\bar{y}_{4}^{2}\bar{y}_{5}^{2}+\bar{y}_{4}\bar{y}_{5}%
\bar{y}_{9}$ \\ 
& $13\bar{y}_{3}^{6}+30\bar{y}_{3}^{2}\bar{y}_{4}^{3}+\bar{y}_{3}\bar{y}_{6}%
\bar{y}_{9}+\bar{y}_{4}^{2}\bar{y}_{5}^{2}+\bar{y}_{4}\bar{y}_{5}\bar{y}_{9}+%
\bar{y}_{6}^{3}$ \\ 
& $-4\bar{y}_{3}^{6}+\bar{y}_{3}^{3}\bar{y}_{4}\bar{y}_{5}-11\bar{y}_{3}^{2}%
\bar{y}_{4}^{3}+\bar{y}_{3}\bar{y}_{6}\bar{y}_{9}+\bar{y}_{5}\bar{y}_{6}\bar{%
y}_{7}$ \\ 
& $13\bar{y}_{3}^{6}+\bar{y}_{3}^{3}\bar{y}_{4}\bar{y}_{5}+28\bar{y}_{3}^{2}%
\bar{y}_{4}^{3}-\bar{y}_{4}^{2}\bar{y}_{5}^{2}+\bar{y}_{4}\bar{y}_{5}\bar{y}%
_{9}$ \\ \hline
\end{tabular}
\end{center}

\noindent \textbf{Step 2. the cohomologies }$H^{\ast }(G/P)$\textbf{. }%
Summarizing the contents of Table 1 we find that

\begin{quote}
$H^{even}(F_{{4}}/P^{s})=\mathbb{Z}[\overline{y}_{{3}},\overline{y}_{{4}},%
\overline{y}_{{6}}]/\left\langle p^{\ast }(h_{{3}}),p^{\ast }(h_{{6}%
}),p^{\ast }(h_{{8}}),p^{\ast }(h_{{12}})\right\rangle $,
\end{quote}

\noindent where

\begin{quote}
$h_{{3}}=2y_{{3}},$ $h_{{6}}=2y_{{6}}+y_{{3}}^{2}$\textsl{, }$h_{{8}}=3y_{{4}%
}^{2}$\textsl{, }$h_{{12}}=y_{{6}}^{2}-y_{{4}}^{3}$,
\end{quote}

\noindent and that $H^{odd}(F_{{4}}/P^{s})$ has the $H^{even}(F_{{4}}/P^{s})$%
-module basis $\{d_{23}\}$. By Theorem 4.3 we obtain the partial presentation

\begin{quote}
$H^{\ast }(F_{{4}}/P)=\mathbb{Z}[\omega _{1},y_{{3}},y_{4},y_{{6}%
}]/\left\langle r_{{3}},r_{6},r_{8},r_{12},r_{12}^{\prime }\right\rangle $,
\end{quote}

\noindent indicating that the inclusion $\left\{ \omega _{1},y_{{3}%
},y_{4},y_{{6}}\right\} \subset H^{\ast }(F_{{4}}/P)$ induces the surjective
ring map $f:\mathbb{Z}[\omega _{1},y_{{3}},y_{4},y_{{6}}]\rightarrow H^{\ast
}(F_{{4}}/P)$. Further, according to Lemma 4.7, computing with the \textsl{%
Null-space }$N(f_{{m}})$ in the order $m=3,6,8$ and $12$ suffices to decide
the generators of the ideal $\ker f$ to yields the Schubert presentation

\begin{enumerate}
\item[(5.11)] $H^{\ast }(F_{4}/P)=\mathbb{Z}[\omega
_{1},y_{3},y_{4},y_{6}]/\left\langle r_{3},r_{6},r_{8},r_{12}\right\rangle $%
, where

$r_{3}=2y_{3}-\omega _{1}^{3}$;$\qquad $

$r_{6}=2y_{6}+y_{3}^{2}-3\omega _{1}^{2}y_{4}$;

$r_{8}=3y_{4}^{2}-\omega _{1}^{2}y_{6}$;

$r_{12}=y_{6}^{2}-y_{4}^{3}$.
\end{enumerate}

Similarly, combining \textsl{the Null-space} with the contents of Tables 2
and 3 one gets the Schubert presentations of the ring $H^{\ast }(G/P)$ for $%
G=$ $E_{6}$ and $E_{7}$ as

\begin{enumerate}
\item[(5.12)] $H^{\ast }(E_{6}/P)=\mathbb{Z}[\omega
_{2},y_{3},y_{4},y_{6}]/\left\langle r_{6},r_{8},r_{9},r_{12}\right\rangle $%
, where

$r_{6}=2y_{6}+y_{3}^{2}-3\omega _{2}^{2}y_{4}+2\omega _{2}^{3}y_{3}-\omega
_{2}^{6};$

$r_{8}=3y_{4}^{2}-6\omega _{2}y_{3}y_{4}+\omega _{2}^{2}y_{6}+5\omega
_{2}^{2}y_{3}^{2}-2\omega _{2}^{5}y_{3};$

$r_{9}=2y_{3}y_{6}-\omega _{2}^{3}y_{6};$

$r_{12}=y_{6}^{2}-y_{4}^{3}.$

\item[(5.13)] $H^{\ast }(E_{7}/P)=\mathbb{Z}[\omega
_{2},y_{3},y_{4},y_{5},y_{6},y_{7},y_{9}]/$\textsl{\ }$\left\langle
r_{j}\right\rangle _{j\in \Lambda }$ where\textsl{\ }$\Lambda
=\{6,8,9,10,12,14,18\}$,

$r_{6}=2y_{6}+y_{3}^{2}+2\omega _{2}y_{5}-3\omega _{2}^{2}y_{4}+2\omega
_{2}^{3}y_{3}-\omega _{2}^{6}$;

$r_{8}=$ \ $3y_{4}^{2}-2y_{3}y_{5}+2\omega _{2}y_{7}-6\omega
_{2}y_{3}y_{4}+\omega _{2}^{2}y_{6}+5\omega _{2}^{2}y_{3}^{2}+2\omega
_{2}^{3}y_{5}-2\omega _{2}^{5}y_{3}$;

$r_{9}=2y_{9}+2y_{4}y_{5}-2y_{3}y_{6}-4\omega _{2}y_{3}y_{5}\ -\omega
_{2}^{2}y_{7}+\omega _{2}^{3}y_{6}\ +2\omega _{2}^{4}y_{5}$;

$r_{10}=y_{5}^{2}-2y_{3}y_{7}+\omega _{2}^{3}y_{7}$;

$r_{12}=y_{6}^{2}+2y_{5}y_{7}-y_{4}^{3}+2y_{3}y_{9}+2y_{3}y_{4}y_{5}+2\omega
_{2}y_{5}y_{6}-6\omega _{2}y_{4}y_{7}+\omega _{2}^{2}y_{5}^{2}$;

$r_{14}=y_{7}^{2}-2y_{5}y_{9}+y_{4}y_{5}^{2}$;

$%
r_{18}=y_{9}^{2}+2y_{5}y_{6}y_{7}-y_{4}y_{7}^{2}-2y_{4}y_{5}y_{9}+2y_{3}y_{5}^{3}-\omega _{2}y_{5}^{2}y_{7} 
$.
\end{enumerate}

\noindent \textbf{Step 3.} \textbf{Computing with the Weyl invariants.\ }In
addition to (5.10)\textbf{\ }the parabolic subgroup $P$ on $G$ specified by
table (5.9) induces also the fibration

\begin{enumerate}
\item[(5.14)] $P/T\overset{i}{\hookrightarrow }G/T\overset{\pi }{\rightarrow 
}G/P$ (i.e. (4.3)),
\end{enumerate}

\noindent where Schubert presentation of the cohomology of the base space $%
G/P$ has been decided by (5.11), (5.12) and (5.13). On the other hand, with

\begin{quote}
$P/T=Sp(3)/T^{3}$, $SU(6)/T^{5}$ or $SU(7)/T^{6}$ for $G=F_{4}$, $E_{6}$ or $%
E_{7}$
\end{quote}

\noindent the cohomology of the fiber space $P/T$ is given by Theorem 5.1 as

\begin{enumerate}
\item[(5.15)] $H^{\ast }(P/T)=\left\{ 
\begin{tabular}{l}
$\frac{\mathbb{Z}[\omega _{2},\omega _{3},\omega _{4}]}{\left\langle
c_{2},c_{4},c_{6}\right\rangle }$ \textsl{if} $G=F_{4}$ \\ 
$\frac{\mathbb{Z}[\omega _{1},\omega _{3},\ldots ,\omega _{n}]}{\left\langle
c_{r},2\leq r\leq n\right\rangle }$ \textsl{if }$G=E_{n}$\textsl{\ with} $%
n=6,7$.%
\end{tabular}%
\right. $
\end{enumerate}

\noindent Thus, Theorem 4.3 is applicable to fashion the ring $H^{\ast
}(G/T) $ in question from the known ones $H^{\ast }(P/T)$ and $H^{\ast
}(G/P) $. To this end we need only to specify a system $\{\rho _{r}\}$
satisfying the constraints (4.5). The invariant theory of Weyl groups serves
this purpose.

Recall that the Weyl group $W$ of $G$ can be identified with the subgroup of 
$Aut(H^{2}(G/T))$ generated by the automorphisms $\sigma _{i}$, $1\leq i\leq
n$, whose action on the Schubert basis $\{\omega _{1},\ldots ,\omega _{n}\}$
of $H^{2}(G/T)$ is given by the Cartan matrix $\left( a_{ij}\right)
_{n\times n}$ of $G$ as

\begin{enumerate}
\item[(5.16)] $\sigma _{{i}}(\omega _{{k}})=\left\{ 
\begin{tabular}{l}
$\omega _{{i}}\text{ if }k\neq i\text{;}$ \\ 
$\omega _{{i}}-\sum\nolimits_{{1\leq j\leq n}}a_{{ij}}\omega _{{j}}\text{ if 
}k=i$%
\end{tabular}%
\right. $.
\end{enumerate}

\noindent Introduce for each $G=F_{4},E_{6}$ and $E_{7}$ the polynomials $%
c_{r}(P)$ in $\omega _{1},_{\cdots },\omega _{n}$ by the formula

\begin{enumerate}
\item[(5.17)] $c_{r}(P):=\left\{ 
\begin{tabular}{l}
$e_{r}(o(\omega _{4},W(P))),1\leq r\leq 4$ if $G=F_{4}$; \\ 
$e_{r}(o(\omega _{n},W(P))$, $1\leq r\leq n$ if $G=E_{n},n=6,7$%
\end{tabular}%
\right. $,
\end{enumerate}

\noindent where $o(\omega ,W(P))\subset H^{2}(G/T)$ denotes the $W(P)$-orbit
through $\omega \in H^{2}(G/T)$, and where $e_{r}(o(\omega ,W(P)))\in $ $%
H^{2r}(G/T)$ is the $r^{th}$ elementary symmetric function on the set $%
o(\omega ,W(P))$. For instance if $G=F_{4},E_{6}$ we have by (5.16) that

\begin{quote}
$o(\omega _{4},W(P))=\left\{ \omega _{_{4}},\omega _{{\footnotesize 3}%
}-\omega _{{\footnotesize 4}},\omega _{{\footnotesize 2}}-\omega _{%
{\footnotesize 3}},\omega _{{\footnotesize 1}}-\omega _{{\footnotesize 2}%
}+\omega _{{\footnotesize 3}},\omega _{{\footnotesize 1}}-\omega _{%
{\footnotesize 3}}+\omega _{{\footnotesize 4}},\omega _{{\footnotesize 1}%
}-\omega _{{\footnotesize 4}}\right\} $,

$o(\omega _{6},W(P))=\left\{ \omega _{{\footnotesize 6}},\omega _{%
{\footnotesize 5}}-\omega _{{\footnotesize 6}},\omega _{{\footnotesize 4}%
}-\omega _{{\footnotesize 5}},\omega _{{\footnotesize 2}}+\omega _{%
{\footnotesize 3}}-\omega _{{\footnotesize 4}},\omega _{{\footnotesize 1}%
}+\omega _{{\footnotesize 2}}-\omega _{{\footnotesize 3}},\omega _{%
{\footnotesize 2}}-\omega _{{\footnotesize 1}}\right\} $.
\end{quote}

\noindent On the other hand, according to Bernstein-Gel'fand-Gel'fand \cite[%
Proposition 5.1]{BGG} the induced map $\pi ^{\ast }$ in (4.3) injects, and
satisfies the relation

\begin{quote}
$\func{Im}\pi ^{\ast }=H^{\ast }(G/T)^{W(P)}=$ $H^{\ast }(G/P)$,
\end{quote}

\noindent implying $c_{r}(P)\in H^{\ast }(G/P)$. Since $c_{r}(P)$ is an
explicit polynomial in the Schubert classes $\omega _{i}$ the \textsl{%
Giambelli polynomials} is functional to express it as a polynomial $g_{r}$
in the special Schubert classes on $H^{\ast }(G/P)$ given by table (5.8):

\begin{enumerate}
\item[(5.18)] 
\begin{tabular}{l||l|l|l}
\hline
$G$ & $F_{4}$ & $E_{6}$ & $E_{7}$ \\ \hline\hline
$g_{2}$ & $4\omega _{1}^{2}$ & $4\omega _{2}^{2}$ & $4\omega _{2}^{2}$ \\ 
$g_{3}$ &  & $2y_{3}+2\omega _{2}^{3}$ & $2y_{3}+2\omega _{2}^{3}$ \\ 
$g_{4}$ & $3y_{4}+2\omega _{1}y_{3}$ & $3y_{4}+\omega _{2}^{4}$ & $%
3y_{4}+\omega _{2}^{4}$ \\ 
$g_{5}$ &  & $3\omega _{2}y_{4}-2\omega _{2}^{2}y_{3}+\omega _{2}^{5}$ & $%
2y_{5}+3\omega _{2}y_{4}-2\omega _{2}^{2}y_{3}+\omega _{2}^{5}$ \\ 
$g_{6}$ & $y_{6}$ & $y_{6}$ & $y_{6}+2\omega _{2}y_{5}$ \\ 
$g_{7}$ &  &  & $y_{7}$ \\ \hline
\end{tabular}%
.
\end{enumerate}

\noindent Up to now we have accumulated sufficient information to show
Theorem 5.4.

\bigskip

\noindent \textbf{Proof of Theorem 5.4.} For each $G=F_{4},E_{6}$ or $E_{7}$
Schubert presentations for the cohomologies of the base $G/P$ and of the
fiber $P/T$ have been determined by (5.11)-(5.13) and (5.15), respectively,
while a system $\{\rho _{r}\}$\ satisfying the relation (4.5) is seen to be $%
\rho _{r}:=c_{r}(P)-g_{r}$. Therefore, Theorem 4.3 is directly applicable to
formulate a presentation of the ring $H^{\ast }(G/T)$. The results can be
further simplified to yield the desired formulae (5.5)-(5.7) by the
following observations:

a) Certain Schubert classes $y_{k}$ from the base space $G/P$ can be
eliminated against appropriate relations of the type $\rho _{k}$, e.g. if $%
G=E_{7}$ the generators $y_{6},y_{7}$ and the relations $\rho _{6},\rho _{7}$
can be excluded by the formulae of $g_{6}$ and $g_{7}$, which implies that $%
y_{6}=c_{6}-2\omega _{2}y_{5}$ and $y_{7}=c_{7}$, respectively;

b) Without altering the ideal, higher degree relations of the type $r_{i}$
may be simplified modulo the lower degree ones by the following fact. For
two ordered sequences $\{f_{i}\}_{1\leq i\leq n}$ and $\{h_{i}\}_{1\leq
i\leq n}$ of a graded polynomial ring with

\begin{quote}
$\deg f_{1}<\cdots <\deg f_{n}$ and $\deg h_{1}<\cdots <\deg h_{n}$
\end{quote}

\noindent write $\{h_{i}\}_{1\leq i\leq n}\thicksim \{f_{i}\}_{1\leq i\leq
n} $ to denote the statements that $\deg h_{i}=\deg f_{i}$ and that $%
(f_{i}-h_{i})\in \left\langle f_{j}\right\rangle _{1\leq j<i}$. Then $%
\{f_{i}\}_{1\leq i\leq n}$ $\sim $ $\{h_{i}\}_{1\leq i\leq n}$ implies that $%
\left\langle h_{1},\ldots ,h_{n}\right\rangle =\left\langle f_{1},\ldots
,f_{n}\right\rangle $.\hfill $\square$

\subsection{A type free characterization of the ring $H^{\ast }(G/T)$}

For an $1$-connected simple Lie group $G$ with rank $n$ denote by $%
D(G)\subset H^{\ast }(G/T)$ the ideal of decomposable elements. Let $h(G)$
be the cardinality of a basis of the quotient group $H^{\ast }(G/T)/D(G)$
and set $m=h(G)-n-1$. The results of Theorems 5.1, 5.2 and 5.4 can be
summarized into one formula, without referring to the types of the group $G$
(see \cite[Theorems 1.2 and 1.3]{DZ3}).

\bigskip

\noindent \textbf{Theorem 5.5.}\textsl{\ For each simple Lie group} $G$ 
\textsl{there exist a set }$\left\{ y_{1},\cdots ,y_{m}\right\} $\textsl{\
of }$m$ \textsl{Schubert classes on }$G/T$\textsl{\ with }$2<\deg
y_{1}<\cdots <\deg y_{m}$\textsl{,} \textsl{so that the inclusion }$\left\{
\omega _{1},\ldots ,\omega _{n},y_{1},\ldots ,y_{m}\right\} \in $\textsl{\ }$%
H^{\ast }(G/T)$ \textsl{induces the Schubert presentation}

\begin{enumerate}
\item[(5.19)] $H^{\ast }(G/T)=\mathbb{Z}[\omega _{1},\ldots ,\omega
_{n},y_{1},\ldots ,y_{m}]/\left\langle e_{i},f_{j},g_{j}\right\rangle
_{1\leq i\leq k;1\leq j\leq m}$\textsl{,}
\end{enumerate}

\noindent \textsl{where}

\begin{quote}
\textsl{i)} $k=n-m$ \textsl{for all }$G\neq E_{8}\QTR{sl}{\ }$\textsl{but} $%
k=n-m+2$ \textsl{for }$G=E_{8}\QTR{sl}{;}$

\textsl{ii) }$e_{i}\in \left\langle \omega _{1},\cdots ,\omega
_{n}\right\rangle $\textsl{, }$1\leq i\leq k$\textsl{;}

\textsl{iii) the pair }$(f_{j},g_{j})$\textsl{\ of polynomials is related to
the Schubert class }$y_{j}$\textsl{\ in the fashion}

$\quad \quad f_{j}$\textsl{\ }$=$\textsl{\ }$p_{j}\cdot y_{j}+\alpha _{j}$%
\textsl{, }$g_{j}=y_{j}^{k_{j}}+\beta _{j}$\textsl{, }$1\leq j\leq m$\textsl{%
, }

\textsl{where }$p_{j}\in \{2,3,5\}$\textsl{\ and }$\alpha _{j},\beta _{j}\in
\left\langle \omega _{1},\cdots ,\omega _{n}\right\rangle $\textsl{;}

\textsl{iv) ignoring the ordering, the sequence} $\left\{ \deg e_{i},\deg
g_{j}\right\} _{1\leq i\leq k;1\leq j\leq m}$ \textsl{of integers agrees
with the degree sequence of the basic Weyl invariants of the group }$G$%
\textsl{\ (over the field} \textsl{of rationals).}\hfill $\square $\noindent
\end{quote}

Concerning assertion iv) we remark that for each simple Lie group $G$ the
degree sequence, as well as explicit formulae, of the basic Weyl invariants%
\textsl{\ }$P_{1},\cdots ,P_{n}$ of $G$ has been determined by Chevalley and
Mehta \cite{Ch0,Me}.

\section{Further remarks on the characteristics}

\textbf{6.1.} Certain parameter spaces of the geometric figures concerned by
Schubert \cite[Chapter 4]{Sc1} may fail to be flag manifolds, but can be
constructed by performing finite number of blow-ups on flag manifolds along
the centers again in flag manifolds, see the examples in Fulton \cite[%
Example 14.7.12]{Fu}, or in \cite{DL} for the construction of the parameter
spaces of the complete conics and quadrics on the $3$-space $\mathbb{P}^{3}$%
. As results the relevant characteristics can be computed from those of flag
manifolds via strict transformations (e.g. \cite[Examples 5.11; 5.12]{DL}).

\textbf{6.2. }As the intersection multiplicities of certain Schubert
varieties on $G/P$, the characteristics $a_{w_{1},\ldots ,w_{k}}^{w}$ are
always non-negative by Van der Waerden \cite{Wa}. Due to the importance of
these numbers in geometry their effective computability (rather than
positivity) had been the top priority in the classical approach \cite%
{Sc1,Sc2,Sc4}, see also Fulton \cite[14.7]{Fu}.

Motivated by the \textsl{Littlewood-Richardson rule} \cite{LR} for the
structure constants of the Grassmannian $G_{n,k}$ a remarkable development
of Schubert calculus has taken place in algebraic combinatorics since
1970's, where the main concern is to find combinatorial descriptions of
characteristics by which the positivity become transparent. This idea has
inspired beautiful results on the enumerations of Yong tableaux, Mondrian
tableaux, Chains in the Bruhat order, and puzzles by A. Buch, W.A. Graham,
I. Coskun, A. Knutson and T.C. Tao \cite{Bu,Co,G,KT1,KT2}, greatly enriched
the classical Schubert calculus.

\textbf{6.3. }According to Van der Waerden \cite{Wa} and A. Weil \cite[p.331]%
{W} Hilbert's 15th problem has been solved satisfactorily. In particular, in
the context of modern intersection theory (e.g. \cite{Fu,GH}) rigorous
treatment of the major enumerative results of Schubert \cite{Sc1} had been
completed independently by many authors (e.g. \cite{A,Kl3,KSX,MP}) \footnote{%
The enumerative results in Schubert \cite{Sc1} were mutually verifiable with
the results of other geometers (e.g. Salmon, Clebsch, Chasles and Zeuthen)
of the same period, hence were already known to be correct at that time.};
granted with the basis theorem the characteristics of flag manifolds can be
evaluated uniformly by the formula (2.2), while the Schubert presentations
of the cohomology rings of flag manifolds have also been available (e.g. 
\cite{B1,Bi,DZ2,DZ3,M}).

However, Schubert calculus remains a vital and powerful tool in constructing
the cohomologies of much broad spaces, such as the homogeneous spaces $G/H$
associated to Lie groups $G$. In contrast to the basis theorem (i.e. Theorem
1.2) the cohomologies of such spaces may be nontrivial in odd degrees, and
may contain torsion elements. Nevertheless, inputting the formula (5.19)
into the Koszul complex

\begin{quote}
$E_{2}^{\ast ,\ast }(G)=H^{\ast }(G/T)\otimes H^{\ast }(T)$
\end{quote}

\noindent associated to the fibration $G\rightarrow G/T$ a unified
construction of the integral cohomology rings of all the $1$-connected
simple Lie groups $G$ has been carried out by Duan and Zhao in \cite{DZ5}.
In addition, the formula (2.2) of the characteristics has been extended to
evaluate the Steenrod operators on the $\func{mod}p$ cohomologies of flag
manifolds \cite{DZ6}, and of the exceptional Lie groups \cite{DZ7}.

\textbf{6.4.} As is of today Schubert calculus has entered the intersection
of several rapidly developing fields of mathematics, and has been
generalized to the studies of other generalized cohomology theories, such as
equivariant, quantum cohomology, K-theory, and cobordism, all of them are
different deformations of the ordinary cohomology \cite{GP}. In this regard
the present paper is by no means a comprehensive survey on the contemporary
Schubert calculus. It illustrates a passage from the Cartan matrices of Lie
groups to the cohomology of homogeneous spaces, where Schubert's
characteristics play a central role..

\bigskip

\noindent \textbf{Acknowledgement.} The authors would like to thank their
referees for valuable suggestions and improvements on the earlier version of
the paper.

\noindent Haibao Duan

Academy of Mathematics and Systems Sciences, Beijing 100190;

School of Mathematical Sciences, University of the Chinese Academy of

Sciences, Beijing 100049,

dhb@math.ac.cn

\bigskip

\noindent Xuezhi Zhao

Department of Mathematics, Capital Normal University, Beijing 100048,

zhaoxve@mail.cnu.edu.cn

\end{document}